\documentclass[11pt]{article}

\usepackage{amsmath, amsthm, amsfonts}
\usepackage[bottom]{footmisc}
\usepackage{geometry}

\newtheorem{defi}{Definition}[section]
\newtheorem{lemma}{Lemma}[section]
\newtheorem{theo}{Theorem}[section]
\newtheorem{prop}{Proposition}[section]
\newtheorem{coro}{Corollary}[section]
\newtheorem{rem}{Remark}[section]
\newtheorem{conj}{Conjecture}[section]

\usepackage{graphicx}
\usepackage{float}
\usepackage{dsfont}
\usepackage{amssymb}
\usepackage{amsmath}
\usepackage{times}
\usepackage{csquotes}

\title{Maximal Spread of Coherent Distributions: a Geometric and Combinatorial Perspective\footnote{University of Warsaw Master's dissertation --- supervised by Dr John M. Noble.}}
\author{Stanis\l{}aw Cichomski\\
\small \\
  \small Faculty of Mathematics, Informatics and Mechanics, \\
  \small University of Warsaw
}

\begin{document}
\maketitle

\begin{abstract}\noindent We discuss some open problems concerning the maximal spread of coherent distributions. We prove a sharp bound on
 $\mathbb{E}|X-Y|^{\alpha}$   for   $(X,Y)$ coherent  and   $\alpha \le 2$,   and establish a novel connection between coherent distributions 
and such combinatorial objects as bipartite graphs, conjugate partitions and Ferrer diagrams.  Our results may turn out to be helpful not only for probabilists, but also for graph theorists, especially for those interested in mathematical chemistry and the study of topological indices.\end{abstract}

\section{Introduction}
\subsection{Background and contributions}
How radically different and contradictory can opinions, stated by two experts or specialists be, while based on distinct  sources of information?
This question,  as in \cite{fundacja},  can be formalised  using the notion of conditional probability. Firstly, both experts must agree upon a basic model of reality, 
which can be understood as accepting common  probability space  $(\Omega, \mathcal{F}, \mathbb{P})$.  Inconsistent sources of information shall  then be identified with different sub  $\sigma$-fields  $\mathcal{G},  \mathcal{H} \subset \mathcal{F}$.  Consequently, opinions involved with judging odds of an event  $A\in \mathcal{F}$,  will be expressed as
random variables  $X, Y,$  defined by
$$X   =   \mathbb{P}(A|\mathcal{G}),$$
$$Y   =   \mathbb{P}(A|\mathcal{H}).$$
Based on \cite{pitman},   we shall refer to such  $(X,Y)$ random vectors as coherent, or alternatively, we might occasionally say that  their joint distribution on  $[0,1]^2$  is coherent.  However, this ambiguity will not lead to any misunderstanding. Note that the characterisation of coherence presented, admits a straightforward extension to vectors of length  $n>2$.  For notational convenience, hereinafter, we write
$$(X, Y)   \in  \mathcal{C},$$
or
$$(X_1, X_2,\cdots,  X_n)  \in  \mathcal{C},$$
whenever we want to indicate, that the vector  $(X,Y)$  or  $(X_1, X_2, \cdots, X_n)$  is coherent. Thus, answering initial question, concerning the maximal spread of coherent opinions, is equivalent to finding (or at least bounding) such quantities, as
$$\sup_{(X,Y)\in \mathcal{C}} \mathbb{E}|X-Y|^{r},$$ 
for  $r\in \mathbb{R}_+$,  or the strongly related quantity
$$\sup_{(X,Y)\in \mathcal{C}}\mathbb{P}(|X-Y|>\delta),$$
for  $\delta \in [0,1]$,   with further variants assuming independence, fixing  $\mathbb{P}(A)$,  e.t.c.. Let us highlight, that this formalism should be
regarded as  taking supremum over all probability spaces  $(\Omega, \mathcal{F}, \mathbb{P})$,  all events  $A\in \mathcal{F}$  and all sub  $\sigma$-fields  $\mathcal{G},  \mathcal{H}  \subset  \mathcal{F}$.  It may seem surprising, but despite the fundamental nature of these problems, they 
have not been studied in depth, at least until lately. In \cite{skew},  expressed differently, it was proved that
\begin{theo} \label{dubins} For all  $n \in \mathbb{Z}_{+}$,  and any  $(X_1, X_2, \cdots, X_n) \in \mathcal{C}$  with  $\mathbb{E}X_1=p$,  we have
$$\mathbb{E}\max_{1\le i \le n} X_i  \ \le \  \frac{p(n-p)}{1+p(n-2)},$$ \end{theo}
\noindent but only recently,  in \cite{pitman},  the following result was established.
\begin{theo} 
$$\sup_{(X, Y)\in \mathcal{C}}\mathbb{E}|X-Y|  =  \frac{1}{2}.$$ \end{theo}
\noindent \textit{Proof}: Fix any  $(X, Y) \in \mathcal{C}$  and let  $p =\mathbb{E}X$.   We use the identity
$$|X-Y|  \ = \  2\cdot \max(X,Y)  -  X  -  Y.$$
Thus, by \ Theorem \ref{dubins}  with  $n=2$,  we have 
$$\mathbb{E}|X-Y|  \ \le  \ 2p(2-p)-2p \  =  \  2p(1-p)  \ \le  \ \frac{1}{2}.$$
To attain the equality, consider  $X'=\mathds{1}_{A}$  and  $Y'=\mathbb{E}{\mathds{1}_{A}}$   for arbitrary  $A\in \mathcal{F}$,  with   $\mathbb{P}(A)=\frac{1}{2}$.  Then 
$$\mathbb{E}|X'-Y'|   =   \mathbb{E}\Big|\mathds{1}_{A}-\frac{1}{2}\Big|   =   \frac{1}{2}. \ \ \  \square$$
It is however doubtful, whether this line of reasoning could be pushed further in order to find
$$\sup_{(X,Y)\in \mathcal{C}} \mathbb{E}|X-Y|^{r},$$ 
for  $r>1$. In fact, one of the main contributions of this thesis is establishing that
$$\sup_{(X,Y)\in \mathcal{C}}\mathbb{E}|X-Y|^{\alpha}   =   2^{-\alpha},$$
for all  $\alpha \in [0,2]$,  which is achieved using only  the  $L^2$-norm  and an elementary geometric framework. Based on this premise, one might suspect that  $2^{-r}$  must turn out to be a  true bound in a general setting. Unfortunately, it is quite easy to construct counterexamples to this hypothesis for  $r>3$.  It seems clear, that progression on this problem for higher exponents, will be associated with establishing some new perspective on theory of coherent opinions. In fact, there are known alternative characterisations of coherent distributions, some of which we shall recall in the next section. Having said that,
let us quote \cite{pitman}  on:
\begin{displayquote}
For reasons we do not understand well, these general characterisations seem to be of little help in establishing the evaluations of $\epsilon(\delta)$ [i.e. $\mathbb{P}(|X-Y|>1-\delta)$] discussed
above, or in settling a number of related problems about coherent distributions [...].
\end{displayquote}
It is our belief, that this is indeed so, because of underlying combinatorial nature of those problems. Let us define
$$\mathcal{C_{I}}  =   \lbrace (X,Y): X,Y \in \mathcal{C} ,  \ X\perp Y\rbrace,$$
as a family  of those coherent distributions, which are additionally independent. Our second important result shows, that for all  $k\in \mathbb{Z}_{+}, 
\ k\ge2,$  we have
$$ \sup_{(X,Y)\in \mathcal{C_{I}}}\mathbb{E}|X-Y|^{k}  =  \sup_{n\in\mathbb{Z}_{+}} \ \sup_{\mathcal{B}(n,n)} \ \frac{1}{n^{2+k}} \sum_{i,j=1}^{n} |\deg(x_{i})-\deg(y_{j})|^{k}, $$
where  $\mathcal{B}(n,n)$  stands for the set of all bipartite graphs with two  $n$  element groups of vertices, i.e.
$$V   =   \{x_{1},\dots,x_{n}\} \ \cup \ \{y_{1},\dots, y_{n}\},$$
and  $\deg(v)$  is a degree of vertex  $v$.  This connection may turn out to be helpful not only for probabilists, but also for graph theorists, especially for those interested in mathematical chemistry and the study of topological indices. For example,  based on \cite{irrt},  for a simple graph  $G$  (i.e. without loops or multiple edges),  one defines its total irregularity measure by
$$\text{irr}_{k,t}(G)   =  \frac{1}{2}\sum_{u,v \in V(G)}|\deg(u)-\deg(v)|^k.$$
Thus, finding supremum of  $\text{irr}_{k,t}$   over all graphs with prescribed number of vertices, seems to be a highly related problem.
To the best of our knowledge, it had not yet arisen broad interest,  with exception of  $k=1,2$.  Note that mentioned graph functionals depend on the choice of  particular graph only through the degree sequence. In this context, we turn our attention to the study of coherent distributions related to conjugate partitions of integers and Ferrer diagrams.  This should not come as  a surprise, since the relation between degree sequences and conjugate partitions is well understood; \ for a comprehensive overview of the topic,   see  \cite{thres}.   Prior to giving formal definitions, let us only mention, that Ferrer diagrams, nowadays attracting growing attention, provide a useful, graphical representations of conjugate partitions. As a slight generalisation of Ferrer diagrams, we define
$$F_{f}   =   \{(u,v)\in[0,1]^2 \ :  \ v < f(u) \} ,$$
where
$$f:[0,1]\rightarrow [0,1],$$
is  any  weakly decreasing step function that takes finitely  many differernt values. We also denote the set of such step functions as   $\mathbb{STEP}.$  Subsequently, for any diagram  $F_f$,  we define a corresponding distribution   $(X_f,Y_f)\in \mathcal{C}_{\mathcal{I}}$,  while ensuring that
$$ \sup_{f \ \in \ \mathbb{STEP}} \ \mathbb{E}|X_{f}-Y_{f}|^{k}   =   \sup_{(X,Y)\in \mathcal{C}_{\mathcal{I}}}\mathbb{E}|X-Y|^{k},$$
for all  $k  \in  \mathbb{Z}_{+}$, \ $k\ge 3$.  Our third and last contribution, is a novel application of those combinatorial ideas to demonstrate that
$$\sup_{f \ \in \ \mathbb{STEP}}\mathbb{P}(|X_f-Y_f|> \delta)   =   2\delta(1-\delta),$$
for  $\delta \in (\frac{1}{2},1]$.  This, at least partially, answers the question raised by   \textit{Burdzy}  and  \textit{Pitman}   in \cite{pitman}, were they have formulated the following conjecture.

\begin{conj}
For  $\delta \in (\frac{1}{2},1]$,  we have
$$\sup_{(X,Y)\in \mathcal{C}_{\mathcal{I}}}\mathbb{P}(|X-Y|\ge \delta)   =   2\delta(1-\delta).$$
\end{conj}
 \noindent As a direct consequence, we also obtain a new upper bound, namely
$$\sup_{(X,Y)\in \mathcal{C}_{\mathcal{I}}}\mathbb{E}|X-Y|^k \  \le  \   2\cdot \frac{k}{(k+1)(k+2)}+  2^{-k}  -  2^{-k-1} \cdot \frac{k(k+3)}{(k+1)(k+2)}.$$

\subsection{Alternative characterisations}
In this section, we provide a short collection of  alternative characterisations of coherent distributions, although these are not referred to elsewhere in the paper.
All of them can be found   in \cite{pitman},  and we refer the interested reader to this excellent resource. 
\begin{prop}\label{P0}Let   $X, Y$   be random variables defined on a probability space  $(\Omega, \mathcal{F}, \mathcal{P})$,  on which one can also define
$$U\sim \mathcal{U}[0,1], \ \ \ \ \ U\perp (X, Y).$$
Then the following conditions are all equivalent:
\begin{enumerate}
\item $(X,Y) \in \mathcal{C}$
\item $0\le X, Y \le 1$,    and for some   $A\in \mathcal{F}$ we have
$$X   =   \mathbb{P}(A|X),$$
$$Y  =  \mathbb{P}(A|Y).$$
\item one can define random variable  $Z$,  with  $0\le Z \le 1$,  such that
$$\mathbb{E}[Zg(X)]   =   \mathbb{E}[Xg(X)],$$
$$\mathbb{E}[Zg(Y)] =  \mathbb{E}[Yg(Y)],$$
for all bounded, measurable functions  $g$  with domain  $[0, 1]$.
\item there exists a measurable function  $\phi:[0,1]^2 \to [0,1]$  such that
$$\mathbb{E}[\phi(X,Y)g(X)]  =  \mathbb{E}[Xg(X)],$$
$$\mathbb{E}[\phi(X,Y)g(Y)]  =  \mathbb{E}[Yg(Y)],$$
for all bounded, measurable functions  $g$  with domain  $[0, 1]$.
\end{enumerate}

\end{prop}

\section{Reduction to bipartite graphs}
\label{chapter: reformulation}

In this chapter we will restate the problem of finding 
 $$\sup_{(X,Y) \in \mathcal{C_{I}}} \mathbb{E}|X-Y|^{2}$$ 
in the language of bipartite graphs. Making use of a graph-theoretic topological index   -  namely the first Zagreb   index  $M_1(G)$  -  we will establish that the solution to the reformulated problem is  $\frac{1}{4}$.

\subsection{Reformulation of the problem}
We start with the definition of independent $\sigma$-fields.
\begin{defi}
Let  $(\Omega, \mathcal{F}, \mathbb{P})$  be a probability space. Let  $\mathcal{F}_1$  and  $\mathcal{F}_2$  be two sub $\sigma$-fields of  $\mathcal{F}$.  Then  $\mathcal{F}_1$ and  $\mathcal{F}_2$   are said to be  independent if for any events  $F_1\in \mathcal{F}_1$  and  $F_2\in \mathcal{F}_2$:
$$\mathbb{P}(F_1\cap F_2)  =  \mathbb{P}(F_1)\cdot \mathbb{P}(F_2).$$
\end{defi}

\noindent The following proposition allows us to perform the discretization.

\begin{prop} \label{Xn} For all   $n\in \mathbb{Z_{+}}$  and any   $(X, Y)\in \mathcal{C}$,  there exists  $(X_{n},Y_{n})\in \mathcal{C}$,  such that  $X_{n}$  and   $Y_{n}$  both take at most  $n$   different values, and 
$$|X-X_{n}|  \le  \frac{1}{n},$$
$$|Y-Y_{n}|  \le  \frac{1}{n}.$$  
Moreover, if   $X\perp Y$,  then we may always choose  $X_{n} \perp Y_{n}$.
\end{prop}
\noindent \textit{Proof}: Fix any  $(X,Y)\in \mathcal{C}$  and assume that it is is defined on the probability space  $(\Omega, \mathcal{H}, \mathcal{P})$.  
Let  $\mathcal{F}$  and  $\mathcal{G}$  be two sub $\sigma$-fields of  $\mathcal{H}$,  such that
$$X  =  \mathbb{P}(A|\mathcal{F})  =  \mathbb{E}(\mathds{1}_{A}|\mathcal{F}),$$
$$Y  =  \mathbb{P}(A|\mathcal{G})  =  \mathbb{E}(\mathds{1}_{A}|\mathcal{G}),$$ 
for some measurable set  $A\in \mathcal{H}$.  Then if  $\sigma_{X}$  is a $\sigma$-field generated by  $X$  and  $\sigma_{Y}$  is a  $\sigma$-field generated by  $Y$,  we have
$$\sigma_{X}\subset \mathcal{F}, \ \ \ \ \ \sigma_{Y}\subset \mathcal{G}.$$
By the tower property we have
$$X  =  \mathbb{E}(X|\sigma_X)  =  \mathbb{E}\Bigg(\mathbb{E}(\mathds{1}_A|\mathcal{F})\Big|\sigma_X\Bigg)  =  \mathbb{E}(\mathds{1}_A|\sigma_X),$$
and similarly
$$Y  =  \mathbb{E}(\mathds{1}_A|\sigma_Y).$$
Let  $\sigma_{X}^{n}$  be the $\sigma$-field geneated by 
$$P_{X}^{n}:=\Bigg\{  \Big\{X\in \Big[0,\frac{1}{n}\Big] \Big\} , \Big\{X\in \Big(\frac{1}{n},\frac{2}{n}\Big]  \Big \} , ...,\Big\{X\in \Big (\frac{n-1}{n},1\Big]  \Big \}       \Bigg\}.$$ 
Set
$$X_{n} =  \mathbb{E}(X|\sigma_{X}^{n})  =  \mathbb{E}\Big(\mathbb{E}(\mathds{1}_{A}  |  \sigma_{X}) \ \Big| \ \sigma_{X}^{n}\Big)  =  \mathbb{E}(\mathds{1}_{A}|\sigma_{X}^{n}),$$
and similarly
$$Y_{n}  =  \mathbb{E}(\mathds{1}_{A}|\sigma_{Y}^{n}),$$
where the last equality in the first line follows from the tower property  and the fact that  $\sigma_{X}^{n}\subset \sigma_{X}$.  Firstly note that, from the above we have  $(X_n,Y_n) \in \mathcal{C}$.  Secondly, since   $P_{X}^{n}$  is $n$-element disjoint partition of   $\Omega$,  $X_{n}$  can take at most  $n$  different values.  Thirdly, by elementary considerations we get  $|X-X_{n}|\le \frac{1}{n}$.  Finally,  independence of   $X_{n}$  and  $Y_{n}$   corresponds to independence of  $\sigma_{X}^{n}$  and  $\sigma_{Y}^{n}$.  But   $\sigma_{X}^{n}\subset \sigma_{X}$,  $\sigma_{Y}^{n}\subset \sigma_{Y}$  imply  $\sigma_{X}^{n} \perp \sigma_{Y}^{n}$  whenever  $\sigma_{X}\perp \sigma_{Y}$. \ \ \  $\square$

\begin{defi}Let   $\mathcal{C}(n)$  be the set of   $(X,Y)\in \mathcal{C}$,  such that  $X$  takes at most  $n$  different values, and  $Y$  takes at most  $n$  different values.  \end{defi}

\begin{defi}Let  $\mathcal{C_{I}}(n)$  be the set of  $(X,Y)\in \mathcal{C_{I}}$,  such that $X$  takes at most  $n$  different values, and  $Y$  takes at most  $n$  different values.  \end{defi}

\begin{prop} \label{newton} We have 
 $$\sup_{(X,Y)\in \mathcal{C}}\mathbb{E}|X-Y|^{2}  =   \sup _{n \in \mathbb{Z}_{+}}\sup_{(X,Y)\in \mathcal{C}(n)}\mathbb{E}|X-Y|^{2}.$$  \end{prop}
\noindent \textit{Proof}: For given  $(X,Y)$  and  $n\in \mathbb{Z}_{+}$ choose  $(X_{n},Y_{n})$  as in previous proposition. Note that
$$ \mathbb{E}|X-Y|^{2}  =  \mathbb{E}|X-X_{n}+X_{n}-Y_{n}+Y_{n}-Y|^{2} \ \le \  \mathbb{E}\Big(|X-X_{n}|+|X_{n}-Y_{n}|+|Y_{n}-Y| \Big)^{2}$$
$$= \ \mathbb{E}\Bigg[ \Big( |X-X_{n}|+|Y-Y_{n}| \Big)^{2}+2|X_{n}-Y_{n}|\Big( |X-X_{n}|+|Y-Y_{n}| \Big) +|X_{n}-Y_{n}|^{2} \Bigg] $$
$$\le \  \frac{4}{n^{2}} +\frac{4}{n}\mathbb{E}|X_{n}-Y_{n}|+ \mathbb{E}|X_{n}-Y_{n}|^{2} \ \le  \ \frac{4}{n^{2}} +\frac{4}{n} + \mathbb{E}|X_{n}-Y_{n}|^{2}.$$  
We can now write
 $$\mathbb{E}|X-Y|^{2}  \ \le  \ \limsup \limits_{n \rightarrow \infty}\Big(\frac{4}{n^{2}} + \frac{4}{n} + \mathbb{E}|X_{n}-Y_{n}|^{2} \Big)$$
 $$ \le  \ \limsup \limits_{n \rightarrow \infty }\Big( \frac{4}{n^2}+\frac{4}{n} \Big)+ \limsup \limits_{n \rightarrow \infty} \mathbb{E}|X_{n}-Y_{n}|^{2},$$
 so  $\mathbb{E}|X-Y|^{2}  \ \le  \ \limsup \limits_{n \rightarrow \infty}\mathbb{E}|X_{n}-Y_{n}|^{2},$  and as a result
$$\sup_{(X,Y)\in \mathcal{C}}\mathbb{E}|X-Y|^{2}  \ \le  \ \sup _{n \in \mathbb{Z}_{+}}\sup_{(X,Y)\in \mathcal{C}(n)}\mathbb{E}|X-Y|^{2}.$$
The inequality in the other direction is clear. \ \ \ $\square$ 
\\ \\ \noindent Repeating the same reasoning with the restriction of independence, gives the following result. 
\begin{coro}  We have
$$\sup_{(X,Y)\in \mathcal{C_{I}}}\mathbb{E}|X-Y|^{2}  =  \sup _{n \in \mathbb{Z}_{+}}\sup_{(X,Y)\in \mathcal{C_{I}}(n)}\mathbb{E}|X-Y|^{2}.$$ \end{coro}

\begin{prop} For every  $n\in \mathbb{Z}_{+}$  we have 
$$\sup_{(X,Y)\in \mathcal{C}(n)}\mathbb{E}|X-Y|^{2}  \ = \ \sup_{A,B}\ \sum_{b_{ij}\neq 0}b_{ij}\Bigg|\frac{\sum_{j}a_{ij}}{\sum_{j}b_{ij}}-\frac{\sum_{i}a_{ij}}{\sum_{i}b_{ij}} \Bigg|^{2},$$
where the supremum is taken over all   $A=(a_{ij}),  B=(b_{ij}) \ \in \ \mathbb{R}^{n \times n}$,  such that  
\begin{equation} \label{w1} \forall_{ij} \ \ 0\le a_{ij} \le b_{ij}, \ \ \ \ \ \ \sum_{ij}b_{ij}=1.   \end{equation} \end{prop}
\noindent \textit{Proof}: Fix  $(X,Y)\in \mathcal{C}(n)$,  and for the time being assume that   $X$  and  $Y$  both take exactly  $n$  different values, namely
$$X(\Omega)  =  \{x_{1}, x_{2}, \cdots, x_{n}\}, \ \ \ \ \ \ Y(\Omega) =  \{ y_{1}, y_{2}, \cdots, y_{n}\}.$$
Now, we can write
$$\sigma_{X}   =   \sigma \Bigg( \Big\{ X^{-1}(x_{1}), X^{-1}(x_{2}), \cdots , X^{-1}(x_{n}) \Big \} \Bigg),$$
$$\sigma_{Y} =   \sigma \Bigg( \Big\{ Y^{-1}(y_{1}), Y^{-1}(y_{2}), \cdots , Y^{-1}(y_{n}) \Big \} \Bigg),$$
which means that  $\sigma$-fields generated by  $X$ and  $Y$  are also generated by given two, disjoint partitions.  
Hence, for all  $1 \le i,j \le n$   and  $\omega \in X^{-1}(x_{i})\cap Y^{-1}(y_{j})$,  we have
$$ X(\omega)  =  \mathbb{E}(\mathds{1}_{A}|\sigma_{X})(\omega)  \ = \  \frac{\mathbb{P}\Big(A\cap \{X=x_{i}\} \Big)}{\mathbb{P}(X=x_{i})}, $$
$$ Y(\omega)  =  \mathbb{E}(\mathds{1}_{A}|\sigma_{Y})(\omega)  \ = \  \frac{\mathbb{P}\Big(A\cap \{Y=y_{j}\} \Big)}{\mathbb{P}(Y=y_{j})}. $$
Thus, setting
$$ \forall_{ij} \ \ \ a_{ij} =  \mathbb{P}\Big(A\cap \{X=x_{i}\} \cap \{Y=y_{j}\}\Big),$$
$$\forall_{ij}  \ \ \  b_{ij}  =  \mathbb{P}\Big(\{X=x_{i}\} \cap \{Y=y_{j}\}\Big),  $$
gives
$$\mathbb{E}|X-Y|^{2} \  =  \ \sum_{b_{ij}\neq 0}b_{ij}\Bigg|\frac{\sum_{j}a_{ij}}{\sum_{j}b_{ij}}-\frac{\sum_{i}a_{ij}}{\sum_{i}b_{ij}} \Bigg|^{2}.$$
Furthermore, for any  $(X,Y)\in \mathcal{C}(n)$   with  $X$  or  $Y$  taking less than  $n$  different values, we can always begin by setting redundant rows or columns of  $A$, $B$   to zero, and  assigning the others as described above. In that way, we have just shown that
$$\sup_{(X,Y)\in \mathcal{C}(n)}\mathbb{E}|X-Y|^{2}  \ \le  \ \sup_{A,B}\ \sum_{b_{ij}\neq 0}b_{ij}\Bigg|\frac{\sum_{j}a_{ij}}{\sum_{j}b_{ij}}-\frac{\sum_{i}a_{ij}}{\sum_{i}b_{ij}} \Bigg|^{2}.$$
To prove the opposite inequality, we start by fixing  $A,B$  such that  (\ref{w1})  holds.
We will give an explicit construction of  $(X', Y') \in \mathcal{C}(n)$,  such that
$$\mathbb{E}|X'-Y'|^{2} \  = \  \sum_{b_{ij}\neq 0}b_{ij}\Bigg|\frac{\sum_{j}a_{ij}}{\sum_{j}b_{ij}}-\frac{\sum_{i}a_{ij}}{\sum_{i}b_{ij}} \Bigg|^{2},$$
defined on the probability space  $([0,1], \mathcal{L}, \lambda ), $  where $\lambda$  is  the lebesgue measure on  $[0,1]$  and  $ \mathcal{L}$  is the $\sigma$-field of  $\lambda$-mesaurable subsets of $[0,1]$.  Start by dividing  $[0,1]$  into a family of  disjoint intervals
 $\{I_{ij}\}_{1\le i, j \le n}$,  such that
$$\forall_{ 1\le i, j \le n} \ \ \ \  \lambda(I_{ij})   =   b_{ij}.$$
For every  $1\le i, j \le n$,   chose a subinterval  $A_{ij}  \subset I_{ij}$,  such that 
$$\forall_{ 1\le i, j \le n} \ \ \ \  \lambda(A_{ij})  =  a_{ij}.$$
Construction of   $\{I_{ij}\}_{1\le i,j \le n}$  and  $\{A_{ij}\}_{1\le i,j \le n}$  is clearly possible by the (\ref{w1})   condition. Set
$$A  =  \bigcup_{1\le i,j \le n}A_{ij},$$
$$\forall_{1\le i\le n } \ \ \ \ G_i  =  \bigcup_{1\le j \le n}I_{ij},$$
$$\forall_{1\le j \le n} \ \ \ \ H_{j}  =  \bigcup_{1\le i \le n}I_{ij}.$$
Thus  $(G_{i})_{1\le i \le n}$  and  $(H_{j})_{1\le j \le n}$  are disjoint partitions of  $[0,1]$,   satisfying 
$$\forall_{ij} \ \ a_{ij}  =  \mathbb{P}(A\cap G_{i} \cap H_{j}),$$
$$\forall_{ij} \ \ b_{ij}  =  \mathbb{P}(G_{i} \cap H_{j}).$$
In this setup, for
$$X'  =  \mathbb{E}\Big[\mathds{1}_{A}\Big|\sigma \Big((G_{i})_{i} \Big)\Big],$$
$$Y'  =  \mathbb{E}\Big[\mathds{1}_{A}\Big|\sigma \Big((H_{j})_{j} \Big)\Big],$$
 we get
$$\mathbb{E}| X'-Y' |^{2}  \ = \ \sum_{b_{ij}\neq 0}b_{ij}\Bigg|\frac{\sum_{j}a_{ij}}{\sum_{j}b_{ij}}-\frac{\sum_{i}a_{ij}}{\sum_{i}b_{ij}} \Bigg|^{2},$$
which completes the proof. \ \ \ $\square$
\\ \\Repeating similar reasoning with the modification of independence, leads to 
\begin{coro} \label{ind} We have
$$\sup_{(X,Y)\in \mathcal{C_{I}}(n)}\mathbb{E}|X-Y|^{2} \ =  \ \sup_{A,B}\ \sum_{b_{ij}\neq 0}b_{ij}\Bigg|\frac{\sum_{j}a_{ij}}{\sum_{j}b_{ij}}-\frac{\sum_{i}a_{ij}}{\sum_{i}b_{ij}} \Bigg|^{2},$$
where supremum is taken over all  $A=(a_{ij}),  B=(b_{ij})  \in  \mathbb{R}^{n \times n}$,  such that  
$$\forall_{ij} \ \ 0\le a_{ij} \le b_{ij},$$
for which there exists  $R=(r_{i}),  C=(c_{j})  \in  \mathbb{R}^{n}$,  satisfying
\begin{equation} \label{w2} \forall_{i} \ \ 0\le r_{i}, \ \ \ \ \  \sum_{i} r_{i}=1,  \end{equation}
$$\forall_{j} \ \ 0\le c_{j}, \ \ \ \ \ \sum_{j} c_{j}=1,$$
$$B  =  RC^{T}.$$ \end{coro}

\begin{defi} Let us define  $\Phi_{n} : [0,1]^{n\times n}\times [0,1]^{n\times n} \longrightarrow \mathbb{R}$  as 
$$\Phi_{n}(A,B) \  =  \ \sum_{ij}\mathds{1}(b_{ij}\neq 0) \cdot b_{ij}\Bigg|\frac{\sum_{j}a_{ij}}{\sum_{j}b_{ij}}-\frac{\sum_{i}a_{ij}}{\sum_{i}b_{ij}} \Bigg|^{2},$$ 
and let   $\mathcal{S}_{n}\subset  [0,1]^{n\times n}\times [0,1]^{n\times n}$  denote the set  of pairs  $(A,B)$  described in  Corollary \ref{ind},  i.e.  satisfying conditions  (\ref{w1})  and  (\ref{w2}).\end{defi}

\begin{prop} \label{compact} The set  $\mathcal{S}_{n}$  is compact. The function  $\Phi_{n}$  is continuous on  $\mathcal{S}_{n}$.\end{prop} 
\noindent \textit{Proof}: It is straightforward to see that  $\mathcal{S}_{n}$  is closed and bounded. To check continuity of  $\Phi_{n}$  it will be enough to verify, that 
$$\forall_{ij} \ \ \ \ \  \phi_{n}^{ij}(A,B) \   :=  \  \mathds{1}(b_{ij}\neq 0)\cdot b_{ij}\Bigg|\frac{\sum_{j}a_{ij}}{\sum_{j}b_{ij}}-\frac{\sum_{i}a_{ij}}{\sum_{i}b_{ij}} \Bigg|^{2},$$
is continuous on  $\mathcal{S}_{n}$.  It is clear, that  $\phi_{n}^{ij}$  is continuous at any  $(A,B)=(a_{kl},b_{kl})_{kl}$  with  $b_{ij}\neq 0$. 
Therfore, let us consider 
$$(A^{(m)},B^{(m)})  \ = \  (a_{kl}^{(m)},b_{kl}^{(m)})_{kl} \ \ \xrightarrow{m\longrightarrow \infty} \ \ (a_{kl},b_{kl})_{kl},$$
with  $b_{ij}=0$.  For   $m$   satisfying   $b_{ij}^{(m)}=0$,  we have   $\phi_{n}^{ij}(A^{(m)},B^{(m)})  =  \phi_{n}^{ij}(A,B)$  $=  0$.  On the other hand,
if  $b_{ij}^{(m)}\neq 0$,  then the given fractions  are well defined, and
$$b_{ij}^{(m)}\cdot \Bigg|\frac{\sum_{j}a_{ij}^{(m)}}{\sum_{j}b_{ij}^{(m)}}-\frac{\sum_{i}a_{ij}^{(m)}}{\sum_{i}b_{ij}^{(m)}} \Bigg|^{2}  \ \le  \ b_{ij}^{(m)}\cdot 1  \ \xrightarrow{m\longrightarrow \infty} b_{ij}=0,$$
which completes the proof. \ \ \ $\square$

\begin{coro} We have
$$\sup_{(X,Y)\in \mathcal{C_{I}}(n)}\mathbb{E}|X-Y|^{2} \  =  \ \sup_{(A,B)\in \mathcal{S}_{n}} \Phi_{n}(A,B).$$\end{coro}

\begin{prop}Without loss of generality, we have
$$ \sup_{\mathcal{S}_{n}}   \Phi_{n}  \ =  \  \sup_{\mathcal{SQ}_{n}} \Phi_{n},$$ 
where  $\mathcal{SQ}_n \subset \mathcal{S}_n $  is the restriction of  $\mathcal{S}_n$  to rationals; namely the subset of those  $(A,B)\in \mathcal{S}_n$,  that satisfy 
$$B  =  RC^T,$$
for some  $R=(r_{i}),  C=(c_{j})  \in  \mathbb{Q}^{n}$,  with
$$\forall_{i} \ \ 0\le r_{i}, \ \ \ \ \  \sum_{i} r_{i}=1, $$
$$\forall_{j} \ \ 0\le c_{j}, \ \ \ \ \ \sum_{j} c_{j}=1.$$\end{prop}
\noindent \textit{Proof}: From  Proposition \ref{compact}  we see, that  $\Phi_{n}$  is uniformly continuous on  $\mathcal{S}_{n}$.  Therefore
 $$\forall  k\in \mathbb{Z}_{+} \ \ \ \exists \delta_{k}>0 \ : \ \ \ \ \forall_{x,y \in \mathcal{S}_{n}} \ \  \ ||x-y||<\delta_{k} \  \implies  \ |\Phi_{n}(x)-\Phi_{n}(y)|<\frac{1}{k}.$$
For given  $(A,B)\in \mathcal{S}_{n}$  and  $k\in \mathbb{Z}_{+}$  choose  $(A_{k},B_{k})\in \mathcal{SQ}_n$  satisfying
$$||(A,B)-(A_{k},B_{k})|| \ < \ \delta_{k}.$$
The set  $\mathcal{SQ}_n$  is clearly dense in  $\mathcal{S}_{n}$  and thus, such  $(A_k,B_k)$  can be found.
Hence, we have
$$\Phi_{n}(A,B) \  <  \ \Phi_{n}(A_{k},B_{k}) + \frac{1}{k} \ \  \implies  \  \ \Phi_{n}(A,B)  \   \le \ \limsup_{k} \ \Phi_{n}(A_{k},B_{k}),$$
and therefore
$$ \sup_{\mathcal{S}_{n}}  \ \Phi_{n} \  \le \   \sup_{\mathcal{SQ}_{n}} \Phi_{n}.$$ 
The inequality in the other direction is clear. \ \ \ $\square$
\\ \\ We will sometimes  omit the subscript and write   $\Phi(x)$   for   $\Phi_{n}(x)$.   By convention, we will also write  $\Phi(A,B)$  for  $A,B \in \mathbb{R}^{m\times n}$  with   $m \neq n$:   we just start by making  $(A,B)$  square matrices first, adding by default zero rows or columns, as needed. 

\begin{defi} For pairs  $A,B\in \mathbb{R}^{m\times n}$  we define the operation  $\Delta_{r}$  of row slicing, as follows: $\forall_ { 1 \le i \le m} \ \forall_{l\in \mathbb{Z}_{+}}$
$$ (A,B)=\Bigg(\begin{bmatrix}
           a_{1} \\ \vdots \\ a_{i-1} \\ a_{i}\\ a_{i+1}\\ \vdots \\ a_{m} \end{bmatrix},
\begin{bmatrix}
           b_{1} \\ \vdots \\  b_{i-1} \\ b_{i}\\ b_{i+1} \\ \vdots \\ b_{m} \end{bmatrix}\Bigg) 
\ \longmapsto \
\Bigg(\begin{bmatrix}
           a_{1} \\ \vdots \\ a_{i-1} \\ a_{i_{1}}\\ \vdots  \\ a_{i_{l}}\\ a_{i+1} \\ \vdots \\ a_{m} \end{bmatrix},
\begin{bmatrix}
           b_{1} \\ \vdots \\ b_{i-1} \\ b_{i_{1}} \\ \vdots \\ b_{i_{l}} \\ b_{i+1} \\ \vdots \\ b_{m} \end{bmatrix}\Bigg)
=\Delta_{r}^{i,l}(A,B),$$
where   $a_{i_{1}}=\dots=a_{i_{l}}=\frac{a_{i}}{l}$  and  $b_{i_{1}}=\dots=b_{i_{l}}=\frac{b_{i}}{l}$. 
We also define the operation  $\Delta_{c}$  of column slicing similarly. \end{defi}

\begin{lemma} \label{cut} Fix  $A,B\in \mathbb{R}^{m\times n}$.  We have
$$ \forall_ { 1 \le i_{0} \le m} \ \ \forall_{l\in \mathbb{Z}_{+}} \ \ \ \ \Phi(A,B)  =  \Phi(\Delta_{r}^{i_{0},l}(A,B)), $$
$$ \forall_ { 1 \le j_{0} \le n} \ \ \forall_{l\in \mathbb{Z}_{+}} \ \ \ \ \Phi(A,B)  =  \Phi(\Delta_{c}^{j_{0},l}(A,B)). $$
\end{lemma}
\noindent \textit{Proof}: We will only prove the first part. Just as before, let us write  $\Phi(A,B)  =  \sum_{ij}\phi^{ij}(A,B)$,  where
$$ \phi^{ij}(A,B)  \  =  \  \mathds{1}(b_{ij}\neq 0)\cdot b_{ij}\Bigg|\frac{\sum_{j}a_{ij}}{\sum_{j}b_{ij}}-\frac{\sum_{i}a_{ij}}{\sum_{i}b_{ij}} \Bigg|^{2}.$$
Start by noting, that row slicing preserves sums of all columns of  $A$  and  $B$.  On the other hand, row slicing can change the sum of a row, only if this particular row was sliced. In this second case, both the corresponding rows of  $A$  and  $B$  have been reduced by the same factor, leaving their proportion unchanged. Therefore we have
$$\forall_{ j, \ i\neq i_{0}} \ \ \ \ \ \phi^{ij}(A,B)  \ = \   \phi^{ij}(\Delta_{r}^{i_{0},l}(A,B)),$$
$$\forall_{ j,  \ i=i_{0}} \ \ \forall_{ 1 \le t \le l} \ \ \ \ \ \phi^{ij}(A,B)\cdot \frac{1}{l}  \ = \  \phi^{i_{t} j}(\Delta_{r}^{i_{0},l}(A,B)),$$
and hence, for $i=i_{0}$
$$\forall_{j} \ \ \ \ \ \phi^{ij}(A,B) \  =  \ \sum_{t=1}^{l} \phi^{i_{t} j}(\Delta_{r}^{i_{0},l}(A,B)).$$
Therefore, summation over the full ranges concludes the proof. \ \ \ $\square$ 
\\ \\Let us  use  $\mathds{1}_{n}$  as notation for $n$-dimensional vector of ones. The following proposition allows us to eliminate the  $b_{ij}$ coefficients.

\begin{prop} We have 
$$ \sup_{n\in \mathbb{Z}_{+}} \ \sup \Big\{ \Phi_{n}(A,B): \ (A,B)\in \mathcal{SQ}_{n} \Big\} $$
$$= \  \sup_{n\in \mathbb{Z}_{+}} \sup \Big\{ \Phi_{n}(A,B): \ (A,B)\in \mathcal{SQ}_n , \ \ B=\frac{1}{n^{2}}\mathds{1}_{n}\mathds{1}_{n}^{T}\Big\}.$$  \end{prop}
\noindent \textit{Proof}: Fix   $(A,B) \in \mathcal{SQ}_{n}$  and   $R=(r_{i}), C=(c_{j}) \in \mathbb{Q}^{n}$,  such that 
$$\forall_{i} \ \ 0\le r_{i}, \ \ \ \ \  \sum_{i} r_{i}=1, $$
$$\forall_{j} \ \ 0\le c_{j}, \ \ \ \ \ \sum_{j} c_{j}=1,$$
$$B  = RC^{T}.$$
Since $\{r_{1},\dots,r_{n}\} \cup \{c_{1},\dots, c_{n}\}$  is a set of rational numbers, there is a common denominator  $D$  and natural numbers 
 $N_{r,1},\dots, N_{r,n}, N_{c,1}, \dots N_{c,n}$,  such that
$$(r_{1},\dots,r_{n}) \  =  \ \Bigg( \frac{N_{r,1}}{D},\dots, \frac{N_{r,n}}{D}\Bigg),$$ 
$$(c_{1},\dots,c_{n}) \  =  \  \Bigg( \frac{N_{c,1}}{D},\dots, \frac{N_{c,n}}{D}\Bigg).$$
Let us now
\\ \\$\bullet$ \ slice every  $i$-th  row of initial  $(A,B)$  matrices exactly  $N_{r,i}$  times,
\\ $\bullet$ \ slice every  $j$-th  column of initial  $(A,B)$  matrices exactly  $N_{c,j}$  times, 
\\ \\ where slicing row or column  $0$  times is to be understood as removing it. 
\\ \\Execution of those operations, leaves us with  $(\tilde{A}, \tilde{B})$,  such that
$$(\tilde{A},\tilde{B})   \in   \mathcal{SQ}_n,$$
$$\tilde{B} =  \frac{\mathds{1}_{\tilde{n}}\mathds{1}^{T}_{\tilde{n}}}{\tilde{n}^{2}},$$
$$\tilde{n} \ = \  \sum_{i=1}^{n}N_{r,i} \ = \   \sum_{j=1}^{n}N_{c,j}.$$
From  Lemma $\ref{cut}$,   it is apparent that   $\Phi(A,B)=\Phi(\tilde{A},\tilde{B})$.  This proves the inequality in one direction. The other direction is clear. \ \ \ $\square$
\\ \\With the analysis so far, we have successfully removed coefficients  $B=(b_{ij})$  from our optimisation problem. Collecting all the pieces together, gives us
\begin{coro}\label{wniosek1} We have
 $$\sup_{(X,Y)\in \mathcal{C_{I}}}\mathbb{E}|X-Y|^{2}  \  =  \ \sup_{n\in \mathbb{Z}_{+}} \ \sup_{A\in [0,1]^{n\times n}} \ \frac{1}{n^{4}}\cdot \sum_{i,j=1}^{n}\Bigg|\sum_{i=1}^{n}a_{ij} - \sum_{j=1}^{n} a_{ij} \Bigg|^{2}.$$ \end{coro}

\begin{defi} Let us define  \ $\Xi_{n}:[0,1]^{n\times n} \longrightarrow \mathbb{R}$  \ as
$$\Xi_{n}(A) \ =  \ \sum_{i,j=1}^{n}\Bigg|\sum_{i=1}^{n}a_{ij} - \sum_{j=1}^{n} a_{ij} \Bigg|^{2}.$$
 \end{defi}
\begin{prop} \label{poprawka}\label{01} For all $n\in \mathbb{Z}_{+}$, we have
 $$\sup_{ A\in[0,1]^{n\times n}}  \Xi_{n}(A) \  =   \sup_{A\in \{0, 1\} ^{n\times n} }   \Xi_{n}(A).$$
 \end{prop}
\noindent \textit{Proof}: Function  $\Xi_{n}$  is continuous on the compact set  $[0,1]^{n\times n}$  and hence it attains its maximum. Let us choose
$$\bar{A}  \ = \  (\bar{a}_{ij})   \ \in \  \underset{[0,1]^{n\times n}}{\operatorname{arg \ max}}  \ \Xi_{n}.$$
For any fixed pair  $(i,j)$  let us set  $\xi_{ij}:[0,1] \rightarrow \mathbb{R}$,
$$\xi_{ij}(a_{ij})  \ =  \ \Xi_{n}(\bar{A} \setminus \bar{a}_{ij}, a_{ij}).$$
The notation means that we use all but one variables of the   $\bar{A}$;   we replace  $\bar{a}_{ij}$  with  $a_{ij}$.  Of course we have
$$\bar{a}_{ij} \ \in  \ \underset{[0,1] }{\operatorname{arg \ max}} \ \xi_{ij}.$$
If   $\bar{a}_{ij}   \not \in  \{0,1\}$,   then  $\xi_{ij}'(\bar{a}_{ij})=0$  and \ $\xi_{ij}''(\bar{a}_{ij})\le 0$.  After some basic calculations, with slight abuse of notation, we get
$$\xi_{ij}'(a_{ij})  \ =  \ 2\cdot \Bigg[n\sum_{j=1}^{n}a_{ij}+n\sum_{i=1}^{n}a_{ij}-2\sum_{i,j=1}^{n}a_{ij}  \Bigg],$$
and
$$\xi_{ij}''(a_{ij}) \ =  \ 2\cdot [n+n-2]  \  \ge  \  \ 0.$$
This proves that, apart from the trivial case  $n=1$,   we cannot have   $\bar{a}_{ij} \in (0,1)$. \ \   $\square$
\\ \\Now, after  Proposition \ref{01},    we can finally explain the connection of our initial problem with bipartite graphs.
\begin{defi} An undirected graph  $G$   is defined as a pair 
 $$G =  (V,E ),$$ 
where  $V$ is a finite set of  vertices and  $E$  is a set of  edges, i.e. unordered pairs of elements of  $V$. \end{defi}
\begin{defi}A simple graph is any undirected graph  $G$,   without loops or multiple edges.
\end{defi}
\begin{defi}A bipartite graph is any simple graph  $G=(V,E)$,  for which  $V$  can be split into two disjoint sets  $V_1$   and  $V_2$,  such that  
each edge   $e\in E$  joins a vertex in  $V_1$  to a vertex in   $V_2$. \end{defi}
\noindent For every  $n\in \mathbb{Z}_{+}$  and  $A\in \{0,1\}^{n\times n}$  consider the graph  $G=(V,E)$  such that  
$$V   =   \{x_{1},\dots,x_{n}\}  \cup  \{y_{1},\dots, y_{n}\},$$
$$ \{x_{1},\dots,x_{n}\}  \cap  \{y_{1},\dots, y_{n}\}  =  \emptyset,$$
$$(x_{i},y_{j})  \in  E   \iff  a_{ij}=1.$$
This leaves us with:

\begin{coro}\label{k=2}$$\sup_{(X,Y)\in \mathcal{C_{I}}}\mathbb{E}|X-Y|^{2}   \ = \  \sup_{n\in \mathbb{Z}_{+}} \ \sup_{\mathcal{B}(n,n)} \ \frac{1}{n^{4}}\cdot \sum_{i,j=1}^{n}|\deg(x_{i})-\deg(y_{j}) |^{2},$$ 
where  $\mathcal{B}(n,n)$  stands for the set of all bipartite graphs with two   $n$   element groups of vertices,
and  $\deg(v)$  is degree of vertex  $v$.\end{coro}
 
\subsection{Solution of the graph problem}
We now show the following:

\begin{theo} \label{graf} For all  $n\in \mathbb{Z}_{+}$,  we have
$$\sup_{\mathcal{B}(n,n)}\sum_{i,j}^{n}|\deg(x_{i})-\deg(y_{j})|^{2}   \ \le  \ \frac{n^{4}}{4}.$$\end{theo}
\noindent Let us start by little simplification
$$\sum_{i,j}^{n}|\deg(x_{i})-\deg(y_{j})|^{2} \ \le \ \frac{n^{4}}{4}  \  \iff  \ n\cdot \Bigg( \sum_{i=1}^{n}\deg^{2}(x_{i})  +\sum_{j=1}^{n} \deg^{2}(y_{j})  \Bigg)$$
$$  \le \ \frac{n^{4}}{4}  +  2\cdot \sum_{i,j=1}^{n}\deg(x_{i})\deg(y_{j}) \  =   \ \frac{n^{4}}{4}  +  2\cdot|E|^{2}, $$
where the last equality follows from  \ $G=(V,E)$ \ being bipartite.

\begin{defi}For any graph  $G=(V,E)$,  we define  first Zagreb index  $M_{1}(G)$,  as
$$M_{1}(G)  =  \sum_{v\in V}\deg^{2}(v).$$ \end{defi}
\noindent A comprehensive overview, of the state of the art of  knowledge of  $M_{1}(G)$, can be found in \cite{zagreb}. In particular we can find there the following

\begin{theo} \label{cudo} Fix $n, e, q \in \mathbb{Z}_{+}$,  $e\le n^{2}$  and let   $e=q\cdot n + r$,  where $0\le r < n$. 
Let  $B^{1}(n,n,e)$  be such a bipartite graph  $G=(V,E)$,  that   $V=X\cup Y$,   $|X|$ $=|Y|=n$,  $|E|=e$,   and 
$q$  vertices from  $Y$  are adjacent to all the vertices in  $X$  and one more vertex from  $Y$  is adjacent to  $r$  vertices
in  $X$.  $B^{1}(n,n,e)$  has its maximum  $M_{1}$  among all  $\mathcal{B}(n,n)$  with  $e$  edges. \end{theo}
\noindent \textit{Proof  of  Theorem  \ref{graf}}:  Take any  $G  =  (V,E)\in \mathcal{B}(n,n)$  with $|E|=e=q\cdot n+r$ as above.
We want to prove that
$$n\cdot M_{1}(G) \ \le \ \frac{n^{4}}{4} + 2(qn+r)^{2}. $$
From  Theorem \ref{cudo}  we can see, that
$$M_{1}(G) \ \le \  (n-r)q^{2}+r(q+1)^{2}+qn^{2}+r^{2},$$
and we simply need to check if
$$n\cdot \Big[(n-r)q^{2}+r(q+1)^{2}+qn^{2}+r^{2}\Big] \ \le \ \frac{n^{4}}{4} + 2(qn+r)^{2}$$
$$\iff \ 0 \ \le \ \frac{n^{4}}{4}-qn^{3}+q^{2}n^{2}+nr(2q-1-r)+2r^{2}$$
$$\iff \ 0 \ \le \Big(\frac{n}{2}-q\Big)^{2}+\Big[\frac{r}{n}(2q-1-r)\Big] +2\Big(\frac{r}{n}\Big)^{2}$$
$$\iff \ 0 \ \le \Big(\frac{n}{2}-q\Big)^{2}+\Big[\frac{r}{n}(2q-n+n-1-r)\Big] +2\Big(\frac{r}{n}\Big)^{2}$$
$$\iff \ 0 \ \le \Big[\Big(q-\frac{n}{2}\Big)+\frac{r}{n} \Big]^{2}+\frac{r}{n}(n-1-r)+\Big(\frac{r}{n}\Big)^{2}.$$
\\The last expression is nonnegative, because  $r+1\le n$   from assumption. \ \ \ $\square$


\section{Solutions for $\mathbb{E}XY$ and $\mathbb{E}|X-Y|^{2}$} 
In this chapter we obtain tight bounds on
$$\sup_{X,Y \in \mathcal{C}(A)}\mathbb{E}XY \ \ \ \ \text{and} \ \ \ \  \sup_{X,Y \in \mathcal{C}(A)}\mathbb{E}|X-Y|^{2},$$
where $\mathcal{C}(A)$ is defined for all $A\in \mathcal{F}$, by 
$$\mathcal{C}(A)  =  \{ \mathbb{E}(\mathds{1}_{A}|\mathcal{G}) \ : \ \mathcal{G}\subset \mathcal{F} \}.$$
Note that, if  $X,Y \in \mathcal{C}(A)$,  then  $(X,Y)$  is clearly coherent. We shall also use
$$\mathcal{C_{I}}(A)  =  \{ (X,Y) \ : \ X,Y\in \mathcal{C}(A), \ X\perp Y \}.$$

\subsection{Two simple bounds on $\mathbb{E}XY$}
To get a better understanding of the definitions, let us start by two exercise-level problems.

\begin{prop} We have
 $$\sup_{(X,Y)\in \mathcal{C_{I}}(A)} \mathbb{E}XY  =  \mathbb{P}(A)^{2}.$$ \end{prop}
\noindent \textit{Proof}: From independence and the tower property of conditional expectation, we get
$$\mathbb{E}XY =  \mathbb{E}X \cdot \mathbb{E}Y  =   \mathbb{P}(A)\cdot \mathbb{P}(A). \ \  \ \square$$

\begin{prop} We have
$$\sup_{(X,Y)\in \mathcal{C}(A)} \mathbb{E}XY  =  \mathbb{P}(A).$$ \end{prop}
\noindent \textit{Proof}: Clearly, for all   $(X,Y)\in \mathcal{C}(A)$  we have  $\mathbb{E}XY  \le  \mathbb{E}X  =  \mathbb{P}(A)$.
Now,  note that
 $$\mathds{1}_{A}   =   \mathbb{E}(\mathds{1}_{A}|\mathcal{F})  \ \in \  \mathcal{C}(A),$$
and hence, putting  $X=Y=\mathds{1}_{A}$,  we get
$$\mathbb{E}XY  =  \mathbb{E}X^{2}  =  \mathbb{E}\mathds{1}_{A}^{2} = \mathbb{E}\mathds{1}_{A}= \mathbb{P}(A). \ \ \ \square$$

\subsection{General bound on $\mathbb{E}|X-Y|^{2}$}
We will start by crystallising the basic geometric intuition in the setting of abstract Hilbert spaces.

\begin{figure}[H]
\centering
\includegraphics[width=70mm]{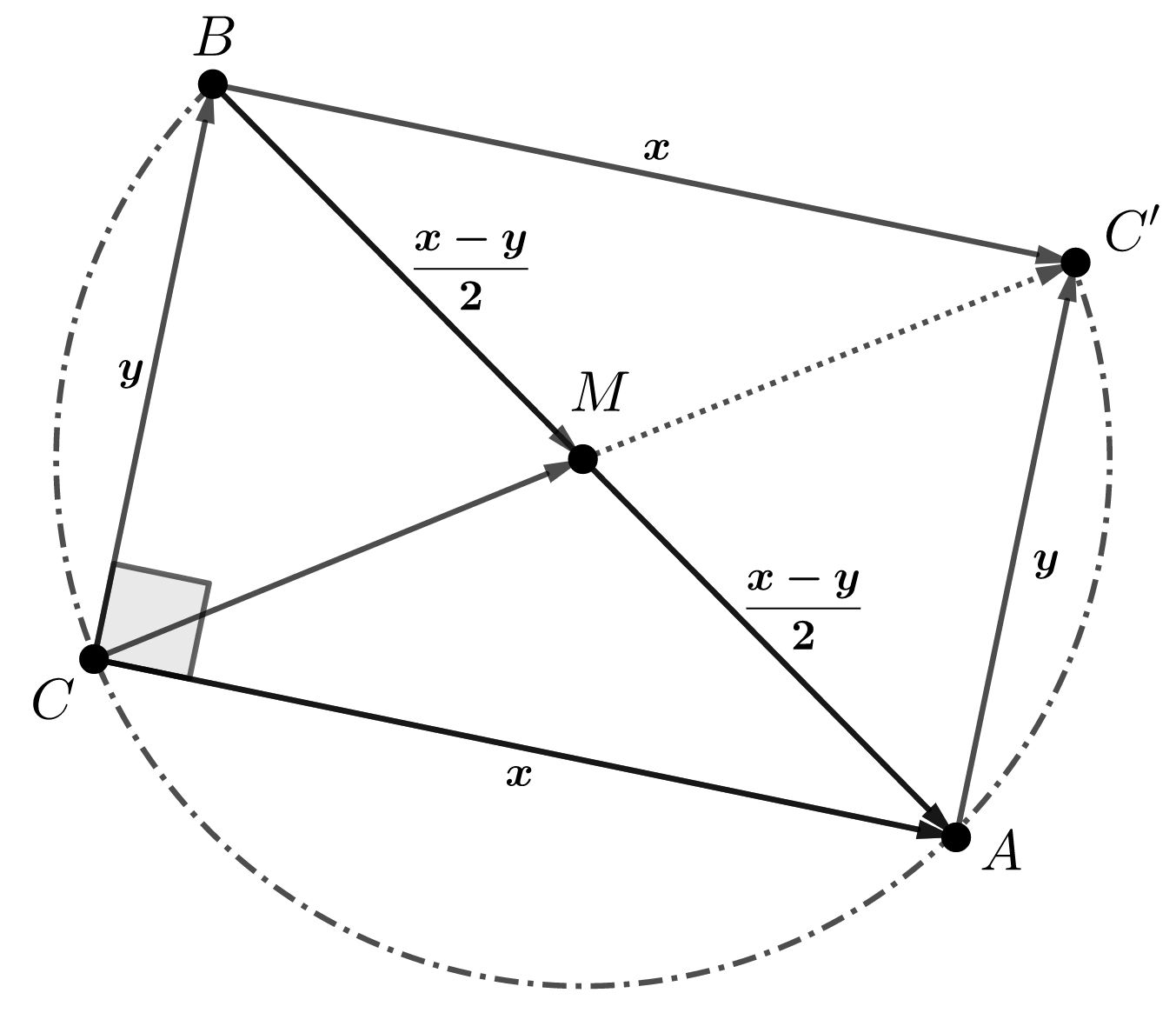}
\caption{ right triangle inscribed in a circle \label{rys1}} \end{figure} 

\begin{lemma} \label{geo} Let  $A,B,C \in L_{2}(\Omega),$   $x=A-C,$   $y=B-C$.  If  $\langle x,y\rangle=0$,   then for  $M=B+\frac{x-y}{2}$   we have  $||B-M||=||C-M||=||A-M||$. \end{lemma}
\noindent \textit{Proof}: From definition, we have $A-B  =  (A-C)-(B-C)  =x-y$ and therefore
$$A-M  =   (A-B)-(M-B)  =  (x-y)-\frac{x-y}{2}   =  \frac{x-y}{2}.$$
This proves that $||B-M||=||A-M||$. Now,  the condition $\langle x,y \rangle = 0$ yields
\begin{equation}||x+y||^{2}  =  ||x||^{2}+2\langle x,y \rangle + ||y||^{2} =  ||x||^{2}-2\langle x,y \rangle + ||y||^{2} \ = \ ||x-y||^{2}. \label{flipflap} \end{equation}
Hence, we get
$$||M-C||=\Big|\Big|(M-B)+(B-C)\Big|\Big|=\Big|\Big|\frac{x-y}{2}+y\Big|\Big|=\Big|\Big|\frac{x+y}{2}\Big|\Big|=||B-M||. \ \ \ \square$$

\begin{prop} \label{alfa} For any  $\alpha \ge 0$,  we have
$$\sup_{(X,Y)\in \mathcal{C}}\mathbb{E}|X-Y|^{\alpha} \ \ge \ \sup_{(X,Y)\in \mathcal{C_{I}}}\mathbb{E}|X-Y|^{\alpha} \ \ge \ 2^{-\alpha}.$$ \end{prop}
\noindent \textit{Proof}: It is enough to set  $X=\mathds{1}_{A}$  and  $Y=\mathbb{E}{\mathds{1}_{A}}$   for arbitrary  $A\in \mathcal{F}$, with   $\mathbb{P}(A)=\frac{1}{2}$.  In such a situation, $Y$  is a constant and therefore  $X\perp Y$.   Then 
$$\mathbb{E}|X-Y|^{\alpha} \ = \ \mathbb{E}\Big|\mathds{1}_{A}-\frac{1}{2}\Big|^{\alpha} \ = \  \frac{1}{2^{\alpha}}. \ \ \  \square$$

\begin{theo} \label{2solve} We have
$$\sup_{(X,Y)\in \mathcal{C}}\mathbb{E}|X-Y|^{2}  =  \frac{1}{4}.$$\end{theo}
\noindent \textit{Proof}:  Fix  the probability space  $(\Omega, \mathcal{F}, \mathbb{P})$  and  $A\in \mathcal{F}$.  We show that 
$$\sup_{X,Y\in \mathcal{C}(A)} \mathbb{E}|X-Y|^{2} \ \le \ \mathbb{P}(A)(1-\mathbb{P}(A))  \ \le \ \Bigg(\frac{\mathbb{P}(A) +(1-\mathbb{P}(A))}{2} \Bigg)^{2}  =  \frac{1}{4}. $$
Start by choosing any two   $\sigma$-fields   $\mathcal{G},\mathcal{H} \subset \mathcal{F}$,  and consider
 $$X = \mathbb{E}(\mathds{1}_{A}|\mathcal{G}), \ \ \ \ \ Y  = \mathbb{E}(\mathds{1}_{A}|\mathcal{H}),$$
$$\mathbb{E}\mathds{1}_{A}  = \mathbb{E}(X|\{\emptyset, \Omega\})  =  \mathbb{E}(Y|\{\emptyset, \Omega\}).  $$

\begin{figure}[H]
\centering
\includegraphics[width=70mm]{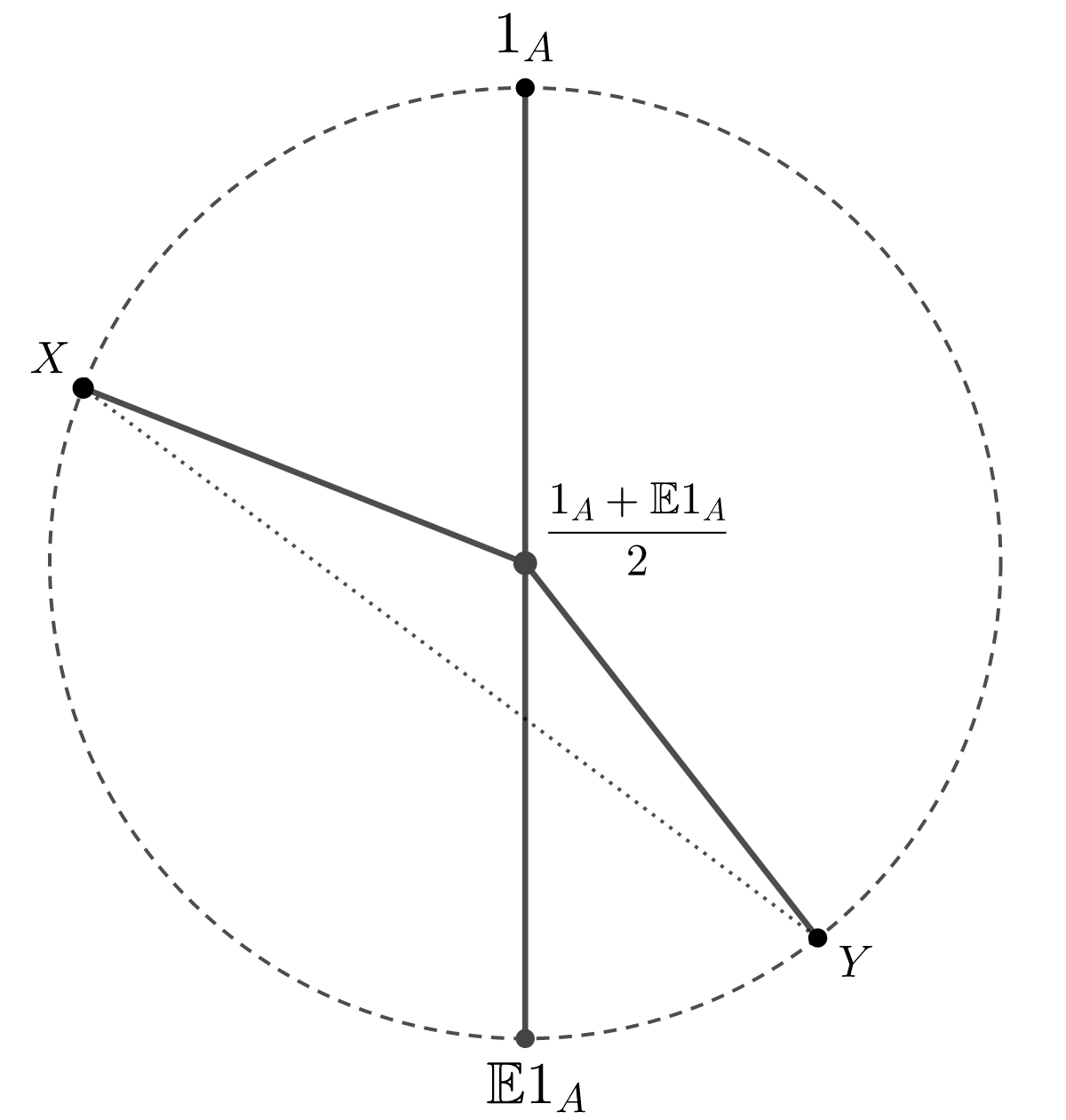}
\caption{ $\mathcal{C}(A)$  and containing it  sphere \label{rys2}} \end{figure}

\noindent
We now check that
$$\Big \langle (\mathds{1}_{A}-X), (X- \mathbb{E}\mathds{1}_{A})\Big\rangle  =  0,$$
namely
$$\Big \langle (\mathds{1}_{A}-X), (X- \mathbb{E}\mathds{1}_{A})\Big\rangle  \ = \  \mathbb{E}\Big[ (\mathds{1}_{A}-X) (X- \mathbb{E}\mathds{1}_{A}) \Big]$$
$$= \  \mathbb{E}\Big[\mathds{1}_{A}X\Big] \ - \ \mathbb{E}\Big[\mathds{1}_A\mathbb{E}(\mathds{1}A)\Big] \ - \ \mathbb{E}\Big[X^2\Big] \ + \ \mathbb{E}\Big[X\mathbb{E}(\mathds{1}_A)\Big]$$
$$= \  \mathbb{E}\Big[\mathds{1}_{A}\mathbb{E}(\mathds{1}_{A}|\mathcal{G})\Big] \ - \ \mathbb{E}\Big[\mathds{1}_A\mathbb{E}(\mathds{1}A)\Big] \ - \ \mathbb{E}\Big[\mathbb{E}(\mathds{1}_{A}|\mathcal{G})^2\Big] \ + \ \mathbb{E}\Big[\mathbb{E}(\mathds{1}_{A}|\mathcal{G})\mathbb{E}(\mathds{1}_A)\Big]$$
$$ =  \ \mathbb{E}\Big[\mathbb{E}(\mathds{1}_{A}|\mathcal{G})^2\Big] \ - \ \mathbb{P}(A)^2 \ - \ \mathbb{E}\Big[\mathbb{E}(\mathds{1}_{A}|\mathcal{G})^2\Big] \ + \ \mathbb{P}(A)^2 \  = \  0.$$  
Similarly
$$\Big\langle (\mathds{1}_{A}-Y), (Y- \mathbb{E}\mathds{1}_{A})\Big\rangle = 0.$$
We have
$$\Bigg|\Bigg|X-\frac{\mathds{1}_{A}+\mathbb{E}\mathds{1}_{A}}{2}\Bigg|\Bigg| \  =  \ \Bigg|\Bigg|\frac{1}{2}(X-\mathds{1}_{A})+\frac{1}{2}(X-\mathbb{E}\mathds{1}_{A})\Bigg|\Bigg|$$
$$ \Bigg|\Bigg|\frac{1}{2}(X-\mathds{1}_{A})-\frac{1}{2}(X-\mathbb{E}\mathds{1}_{A})\Bigg|\Bigg| \  =  \ \frac{1}{2}\cdot||\mathds{1}_{A}-\mathbb{E}\mathds{1}_{A}||,$$
were we have flipped the sign by observation  (\ref{flipflap}).  Similarly
$$ \Bigg|\Bigg|Y-\frac{\mathds{1}_{A}+\mathbb{E}\mathds{1}_{A}}{2}\Bigg|\Bigg| \ =  \ \frac{1}{2}\cdot||\mathds{1}_{A}-\mathbb{E}\mathds{1}_{A}||.$$
Applying the triangle inequality, we get
$$||X-Y||  \ \le  \ \Bigg|\Bigg|X-\frac{\mathds{1}_{A}+\mathbb{E}\mathds{1}_{A}}{2}\Bigg|\Bigg| \ + \ \Bigg|\Bigg|\frac{\mathds{1}_{A}+\mathbb{E}\mathds{1}_{A}}{2}-Y\Bigg|\Bigg| $$
 $$= \  \frac{1}{2}\cdot||\mathds{1}_{A}-\mathbb{E}\mathds{1}_{A}|| + \frac{1}{2}\cdot||\mathds{1}_{A}-\mathbb{E}\mathds{1}_{A}||  \ =  \ ||\mathds{1}_{A}-\mathbb{E}\mathds{1}_{A}||,$$\
resulting in
$$\mathbb{E}|X-Y|^{2} \ \le \ \mathbb{E}|\mathds{1}_{A}-\mathbb{E}\mathds{1}_{A}|^{2} \ = \ (1-\mathbb{P}(A))^{2}\cdot \mathbb{P}(A)+\mathbb{P}(A)^{2}\cdot (1-\mathbb{P}(A)) $$
$$= \  \mathbb{P}(A)(1-\mathbb{P}(A)),$$
which completes the proof. \ \ \ $\square$
\\ \\The following corollary is immediate from the analysis above.

\begin{coro}  \label{sferka} For fixed probability space  $(\Omega, \mathcal{F}, \mathbb{P})$  and an event  $A\in \mathcal{F}$,  we have
$$ \mathcal{C}(A) \ \subset \ \Bigg \{ X\in L_2(\Omega) \ : \  \Bigg|\Bigg|X-\frac{\mathds{1}_{A}+\mathbb{E}\mathds{1}_{A}}{2}\Bigg|\Bigg|
= \frac{\sqrt{\mathbb{P}(A)(1-\mathbb{P}(A))}}{2} \Bigg \}.$$ \end{coro}

\begin{coro} For all  $\alpha \in [0,2]$,  we have
$$\sup_{(X,Y)\in \mathcal{C}}\mathbb{E}|X-Y|^{\alpha} \ = \ \sup_{(X,Y)\in \mathcal{C_{I}}}\mathbb{E}|X-Y|^{\alpha} \ =\ 2^{-\alpha}.$$ \end{coro}
\textit{Proof}: Thanks to \ Proposition $\ref{alfa}$, \ we only need to verify that
$$\sup_{(X,Y)\in \mathcal{C}}\mathbb{E}|X-Y|^{\alpha} \ \le \ 2^{-\alpha}.$$
Clearly \ $\frac{\alpha}{2} \in [0,1]$,\ and thus \ $f(x)=x^{\frac{\alpha}{2}}$ \ is concave on $\mathbb{R}_{+}$. By Jensen inequality, we get
$$\mathbb{E}|X-Y|^{\alpha} \ \le \ \Big( \mathbb{E}|X-Y|^{2} \Big)^{\frac{\alpha}{2}} \ \le \ \Bigg(\frac{1}{4}\Bigg)^{\frac{\alpha}{2}}\ = \ 2^{-\alpha},$$
for all $(X,Y) \in \mathcal{C}$. \ \ \ $\square$

\subsection{Geometry of the multivariate case}
In this section we will obtain an upper bound on 
$$ \sup_{(X_1,\cdots, X_n)  \ \in \ \mathcal{C}(A)} \  \frac{1}{2}\sum_{\substack{i,j=1 \\ i\not=j}}^{n}\mathbb{E}|X_i-X_j|^{2}.$$
By  Theorem  \ref{2solve},  we could simply write
$$ \sup_{(X_1,\cdots, X_n)  \ \in \ \mathcal{C}(A)} \  \frac{1}{2} \sum_{\substack{i,j=1 \\ i\not=j}}^{n}\mathbb{E}|X_i-X_j|^{2} \   \le  \ \frac{n(n-1)}{2}\cdot \mathbb{P}(A)(1-\mathbb{P}(A)),$$
but it turns out  that by using geometric tools,   we can improve it  by a factor of roughly  $2$:
$$ \sup_{(X_1,\cdots, X_n)  \ \in \ \mathcal{C}(A)}  \ \frac{1}{2} \sum_{\substack{i,j=1 \\ i\not=j}}^{n}\mathbb{E}|X_i-X_j|^{2} \   \le \  \frac{n^2}{4} \cdot \mathbb{P}(A)(1-\mathbb{P}(A)).$$
We start with the observation, that enables us to work in the much more intuitive space  $\mathbb{R}^{n}$,  rather than the abstract Hilbert space  $L_{2}(\Omega)$.

\begin{prop} \label{isometria} For any   $(X_1,X_2, \cdots, X_n)   \in  \mathcal{C}(A),$   there are   $x_1,x_2,\cdots, x_n \in \mathbb{R}^{n}$, 
such that 
 \\ \\ $a)$ \ for all    $ 1 \le i,j \le n$   we have $\mathbb{E}|X_i-X_j|^2=||x_i-x_j||^{2}$,
\\ \\ $b)$ \ the set  $\{x_1,x_2,\cdots, x_n\}$   lies on a sphere with  radius   $\frac{ \sqrt{\mathbb{P}(A)(1-\mathbb{P}(A))}}{2}$.
 \end{prop}
\noindent \textit{Proof}: The random variables  $\Big \{X_1, X_2,\cdots, X_n, \frac{\mathds{1}_{A}+\mathbb{E}\mathds{1}_{A}}{2}\Big \}$    are at most  $n+1$  different points in the Hilbert space  $L_{2}(\Omega).$  Therefore they must lie on an  $n$-dimensional affine subspace $H$. Since it is finite dimensional,  $H$  is isometric to the euclidean space  $\mathbb{R}^{n}$.  Let  $x_1, x_2, \cdots, x_n $  be the respectiv images of  $X_1, X_2, \cdots, X_n$  under this isometry. Point  $a)$  then follows automatically   and   point  $b)$  is a direct consequence of  Corollary \ref{sferka}. \ \ \ $\square$
\\ \\By  Proposition \ref{isometria},  we get the following geometric restriction 

\begin{coro}\label{gwoździe} We have  $$ \sup_{(X_1,\cdots, X_n)  \ \in \ \mathcal{C}(A)}  \ \frac{1}{2} \sum_{\substack{i,j=1 \\ i\not=j}}^{n}\mathbb{E}|X_i-X_j|^{2}  \  \le  \ \sup _{x_1,\cdots, x_n \ \in  \ \mathcal{S}} \ \frac{1}{2}\sum_{\substack{i,j=1 \\ i\not=j}}^{n}||x_i-x_j||^{2},$$
where   $\mathcal{S}\subset \mathbb{R}^n$ is a sphere with  a radius   $\frac{ \sqrt{\mathbb{P}(A)(1-\mathbb{P}(A))}}{2}$. \end{coro}

\begin{defi} Let  $\mathcal{M}=\{x_1,\cdots,x_m\}$  be a fininte multiset of points  in  $\mathbb{R}^{n}$.   We will say, that  $\bar{x}  \in   \mathbb{R}^{n}$ is  a  mass centre   of  $\mathcal{M}$,   if
$$\bar{x}  =  \frac{1}{m}\sum_{i=1}^{m}x_i.$$ \end{defi}
\noindent We now recall  Definition \ref{chord} and  Theorem \ref{npoints},    which can be found in  \cite{kwadraty}.

\begin{defi}\label{chord} Let  $\mathcal{M}$  be a multiset of points on an  $(n-1)$-sphere  in  $\mathbb{R}^{n}$.  We define a chord of 
$\mathcal{M}$  to be a line segment whose endpoints belong to  $\mathcal{M}$. \end{defi}

\begin{theo}\label{npoints} Let  $\mathcal{M}$ be a multiset of   $m$  points on a unit  $(n-1)$-sphere,   and let
 $\mathcal{C}$  be the multiset of the lengths of all the chords between them. Then
$$\sum_{c\in \mathcal{C}}c^2  =   m^2(1-d^2),$$ 
where  $d$  is the distance between the  mass centre  of  $\mathcal{M}$  and the centre of the unit
 $(n-1)$-sphere.\end{theo}
\noindent Let us emphasise an important feature of Theorem \ref{npoints}  with the following remark.
\begin{rem}\label{kilof}  The  sum  $\sum_{c\in \mathcal{C}}c^2$  depends on the configuration of  the
 $\{x_1,\cdots,x_m\}$ $=\mathcal{M}$  only through the number of points  $m$  and mass centre  $\bar{x}$.  It does not depend on the affine dimmension of  $\mathcal{M}$. \end{rem}

\begin{theo} \label{multi} We  have
$$ \sup_{(X_1,\cdots, X_n)  \ \in \ \mathcal{C}(A)} \  \frac{1}{2} \sum_{\substack{i,j=1 \\ i\not=j}}^{n}\mathbb{E}|X_i-X_j|^{2} \   \le  \ \frac{n^2}{4} \cdot \mathbb{P}(A)(1-\mathbb{P}(A)).$$ \end{theo}
\noindent \textit{Proof}: By  Corollary \ref{gwoździe},  it suffices to show that
$$\frac{1}{2} \sum_{\substack{i,j=1 \\ i\not=j}}^{n}||x_i-x_j||^{2} \  =  \ \frac{1}{2}\sum_{i,j=1}^{n}||x_i-x_j||^{2}  \ \le  \ \frac{n^2}{4} \cdot \mathbb{P}(A)(1-\mathbb{P}(A)),$$
for all  $x_1,\cdots,x_n \in \mathbb{R}^{n}$  lying on a sphere with a radius  $r=\frac{ \sqrt{\mathbb{P}(A)(1-\mathbb{P}(A))}}{2}$,
which is an immediate consequence of  Theorem \ref{npoints} and scaling by a factor  $r^2$.   $\square$ 
\\ \\It is not clear wether the inequlaity in  Theorem \ref{multi}  can be replaced by an equality sign for all  $n$,  it is however straightforward to attain equality for even   $n=2k$;  for   $i \in \{1,2,\cdots, 2k\}$,  set
$$X_i=\left\{ \begin{array}{ll}
\mathds{1}_A & \text{for} \ \ \ 2 \ | \ i\\
 \mathbb{P}(A) & \text{for} \ \ \ 2 \nmid i,
\end{array} \right. $$
which can be thought of as placing  an equal number of points on each side of the diameter. 

\section{Relationship with Ferrer diagrams}
In this chapter, exploiting the graphical representation introduced earlier, we will establish a connection between 
$$\sup_{(X,Y) \in \mathcal{C_{I}}} \mathbb{E}|X-Y|^{k},$$
for  $k\in \mathbb{Z}_{+}$, $k\ge3$,   and well studied conjugate partitions of integers.

\subsection{Reduction to bipartite graphs:  $k\ge 3$}
\label{section: fromktobi}

\begin{theo} \label{Bk} For all  $k\in \mathbb{Z}_{+},  k\ge3,$  we have
$$ \sup_{(X,Y)\in \mathcal{C_{I}}}\mathbb{E}|X-Y|^{k} \ = \ \sup_{n\in\mathbb{Z}_{+}} \ \sup_{\mathcal{B}(n,n)} \ \frac{1}{n^{2+k}} \sum_{i,j=1}^{n} |\deg(x_{i})-\deg(y_{j})|^{k}. $$\end{theo}
\noindent To prove this result it is enough to reconsider  argumentation presented in  Chapter 1. Only  Proposition \ref{newton}  and  Proposition \ref{poprawka}  used the assumption  $k=2$  explicitly. In the rest of this section we will show that analogous statements hold for any $k$.

\begin{prop}   We have   $$\sup_{(X,Y)\in \mathcal{C_{I}}}\mathbb{E}|X-Y|^{k}  \ =  \ \sup _{n \in \mathbb{Z}_{+}}\ \sup_{(X,Y)\in \mathcal{C_{I}}(n)}\ \mathbb{E}|X-Y|^{k}.$$  \end{prop} 
\noindent \textit{Proof}: For given $(X,Y)$ and  $n\in \mathbb{Z}_{+}$ choose $(X_{n},Y_{n})$ as in  Proposition \ref{Xn}. Note that
$$ \mathbb{E}|X-Y|^{k}=\mathbb{E}|X-X_{n}+X_{n}-Y_{n}+Y_{n}-Y|^{k} \le \mathbb{E}\Big(|X-X_{n}|+|X_{n}-Y_{n}|+|Y_{n}-Y| \Big)^{k}$$
$$= \  \sum_{j=0}^{k} \binom{k}{j} \mathbb{E}\Bigg[\Big(|X-X_{n}|+|Y-Y_{n}| \Big)^{k-j}|X_{n}-Y_{n}|^{j}\Bigg]$$
$$\le  \ \sum_{j=0}^{k-1} \binom{k}{j}\mathbb{E}\Bigg[\Bigg(\frac{1}{n}+\frac{1}{n}\Bigg)^{k-j}\cdot |X_{n}-Y_{n}|^{j} \Bigg]  +  \mathbb{E}|X_{n}-Y_{n}|^{k}.$$
$$\le \   \sum_{j=0}^{k-1} \binom{k}{j}\cdot \Bigg(\frac{2}{n}\Bigg)^{k-j}   +  \mathbb{E}|X_{n}-Y_{n}|^{k}.$$
We can now write
 $$\mathbb{E}|X-Y|^{k} \ \le \ \limsup \limits_{n \rightarrow \infty}\Bigg(\sum_{j=0}^{k-1}\binom{k}{j}\cdot  \Bigg(\frac{2}{n}\Bigg)^{k-j} + \mathbb{E}|X_{n}-Y_{n}|^{k} \Bigg) $$
$$ \le \ \limsup \limits_{n \rightarrow \infty }\Bigg[ \ \sum_{j=0}^{k-1}\binom{k}{j}\cdot  \Bigg(\frac{2}{n}\Bigg)^{k-j}  \ \Bigg]  +   \limsup \limits_{n \rightarrow \infty} \mathbb{E}|X_{n}-Y_{n}|^{k},$$
so  $\mathbb{E}|X-Y|^{k} \ \le \ \limsup \limits_{n \rightarrow \infty}\mathbb{E}|X_{n}-Y_{n}|^{k}$  and as a result
$$\sup_{(X,Y)\in \mathcal{C}}\mathbb{E}|X-Y|^{k} \ \le \ \sup _{n \in \mathbb{Z}_{+}}\sup_{(X,Y)\in \mathcal{C}(n)}\mathbb{E}|X-Y|^{k}.$$
Inequality in the other direction is clear. \ \ \  $\square$

\begin{defi} Let us define  $\Xi_{n}^{k}:[0,1]^{n\times n} \longrightarrow \mathbb{R}$  as
$$\Xi_{n}^{k}(A)  =  \sum_{i,j=1}^{n}\Bigg|\sum_{i=1}^{n}a_{ij} - \sum_{j=1}^{n} a_{ij} \Bigg|^{k}.$$
For  $i,j \in \{1,2,\cdots,n\}$  we  introduce the abbreviation  
$$A_{i\bullet}  =  \sum_{j=1}^{n}a_{ij},$$
$$A_{\bullet j}  =  \sum_{i=1}^{n}a_{ij}.$$
We can now write
$$\Xi_{n}^{k}(A)  =  \sum_{i,j=1}^{n}|A_{i\bullet}-A_{\bullet j}|^{k}.$$
 \end{defi}

\begin{lemma}\label{|x|} For any   $x\in \mathbb{R}$   and  $k\ge 3$  we have 
$$\frac{\partial }{\partial x}|x|^{k}  =  k|x|^{k-2}\cdot x,$$
$$\frac{\partial}{\partial x}(|x|^{k-2}\cdot x)  =  (k-1)|x|^{k-2}.$$
 \end{lemma}
\noindent \textit{Proof}:
$$
\frac{\partial}{\partial x} |x|^{k}  \  \ =  \ \
\begin{cases}
\frac{\partial}{\partial x} x^{k} &\text{:} \ \ x\ge 0\\
\frac{\partial}{\partial x}(-1)^{k}x^{k} &\text{:}\ \ x\le 0
\end{cases} \ \ = \ \
\begin{cases}
kx^{k-1} &\text{:} \ \ x\ge 0\\
(-1)^{k}kx^{k-1} &\text{:}\ \ x\le 0
\end{cases} \ \  = \ \ k|x|^{k-2}\cdot x,
$$
$$
\frac{\partial}{\partial x} (|x|^{k-2}\cdot x) \ \ = \ \
\begin{cases}
\frac{\partial}{\partial x} x^{k-1} &\text{:} \ \ x\ge 0\\
\frac{\partial}{\partial x}(-1)^{k-2}x^{k-1} &\text{:}\ \ x\le 0
\end{cases} \ \ = \ \
\begin{cases}
(k-1)x^{k-2} &\text{:} \ \ x\ge 0\\
(-1)^{k-2}(k-1)x^{k-2} &\text{:}\ \ x\le 0
\end{cases} \ \ =\ \ (k-1)|x|^{k-2},
$$
which ends the proof. \ \ \ $\square$

\begin{prop} \label{01} For all $n\in \mathbb{Z}_{+}$, we have
 $$\sup_{ [0,1]^{n\times n}}  \Xi_{n}^{k}   \ = \   \sup_{\{0, 1\} ^{n\times n} }   \Xi_{n}^{k}.$$
 \end{prop}
\noindent \textit{Proof}:  The function  $\Xi_{n}^{k}$  is continuous on the compact set  $[0,1]^{n\times n}$, and hence it attains its maximum. Let us choose
$$\bar{A}  =  (\bar{a}_{ij})  \ \in  \ \underset{[0,1]^{n\times n}}{\operatorname{arg \ max}} \ \Xi_{n}^{k}.$$
For any fixed pair  $(i,j)$,  let us put  $\xi_{ij}:[0,1] \rightarrow \mathbb{R}$,
$$\xi_{ij}(a_{ij})  \ = \  \Xi_{n}^{k}(\bar{A} \setminus \bar{a}_{ij}, a_{ij}),$$
meaning that we have all but one of the variables as in  $\bar{A}$,   we replace  $\bar{a}_{ij}$  with $a_{ij}$.  With a slight abuse of notation, we have
$$\frac{\partial}{\partial a_{ij}}\xi_{ij}(a_{ij})  =  \frac{\partial}{\partial a_{ij}}\Big(\sum_{p\not=j}|A_{i\bullet}-A_{\bullet p}|^{k} \ + \ \sum_{p\not=i}|A_{p\bullet}-A_{\bullet j}|^{k}  +  |A_{i\bullet}-A_{\bullet j}|^{k}\Big).$$
Note that   $a_{ij}$  cancels out in   $|A_{i\bullet}-A_{\bullet j}|$. By Lemma \ref{|x|},   we now have
$$\frac{\partial}{\partial a_{ij}}\xi_{ij}(a_{ij}) \ = \ k\Bigg[\sum_{p\not=j}|A_{i\bullet}-A_{\bullet p}|^{k-2}(A_{i\bullet} -A_{\bullet p}) \ - \ \sum_{p\not=i}|A_{p\bullet}-A_{\bullet j}|^{k-2}(A_{p\bullet}-A_{\bullet j}) \Bigg],$$
$$\frac{\partial ^{2}}{\partial a_{ij}^2}\xi_{ij}(a_{ij}) \ = \ k(k-1)\Bigg[ \sum_{p\not=j}|A_{i\bullet}-A_{\bullet p}|^{k-2}\  + \ \sum_{p\not=i}|A_{p\bullet}-A_{\bullet j}|^{k-2} \Bigg] \ \ge \ 0,$$
and hence  $\xi_{ij}$ is a convex function. Since maximum of a convex function on compact, convex set is attained on the boundary, we can without loss of generality assume, that  $\bar{a}_{ij}\in \{0,1\}$  for all  $i,j$. \ \ \ $\square$

\subsection{Degree sequences and majorization}
For  $n\in \mathbb{Z}_{+}$  and given two integer sequences  $a=(a_{i})_{i=1}^{n}$, $b=(b_{j})_{j=1}^{n}$,  with
$$n  \ \ge  \  a_{1}, \  a_{2},  \ \cdots , \ a_{n} \  \ge  \ 0,$$   
$$n  \ \ge \  b_{1}, \  b_{2}, \   \cdots , \  b_{n} \  \ge \  0,$$   
we might wonder if there exists a bipartite graph  $G\in \mathcal{B}(n,n)$  with   degree sequences  $a$   and  $b$  in each part, respectively.
We shall call such  $(a, b)$   pairs bigraphic. This question can be answered by famous   Theorem \ref{Ross-Geller} \textit{(Gale-Ryser)},  see for example  \cite{g-r}  or  \cite{thres}.

\begin{defi} \label{cp} For   $n\in \mathbb{Z}_{+}$  and any integer sequnece   $b=(b_{i})_{i=1}^{n}$,  with
$$n  \ \ge  \  b_{1}, \  b_{2},  \ \cdots , \ b_{n} \  \ge  \ 0,$$   
we define its conjugate partition $b^{*}=(b^{*}_{i})_{i=1}^{n}$,   by
$$b^{*}_{k}  =  |\{i:b_{i}\ge k\}|,$$
for all   $k\in \{1,2,\cdots,n\}$.
\end{defi}

\begin{figure}[H]
\centering
\includegraphics[width=55mm]{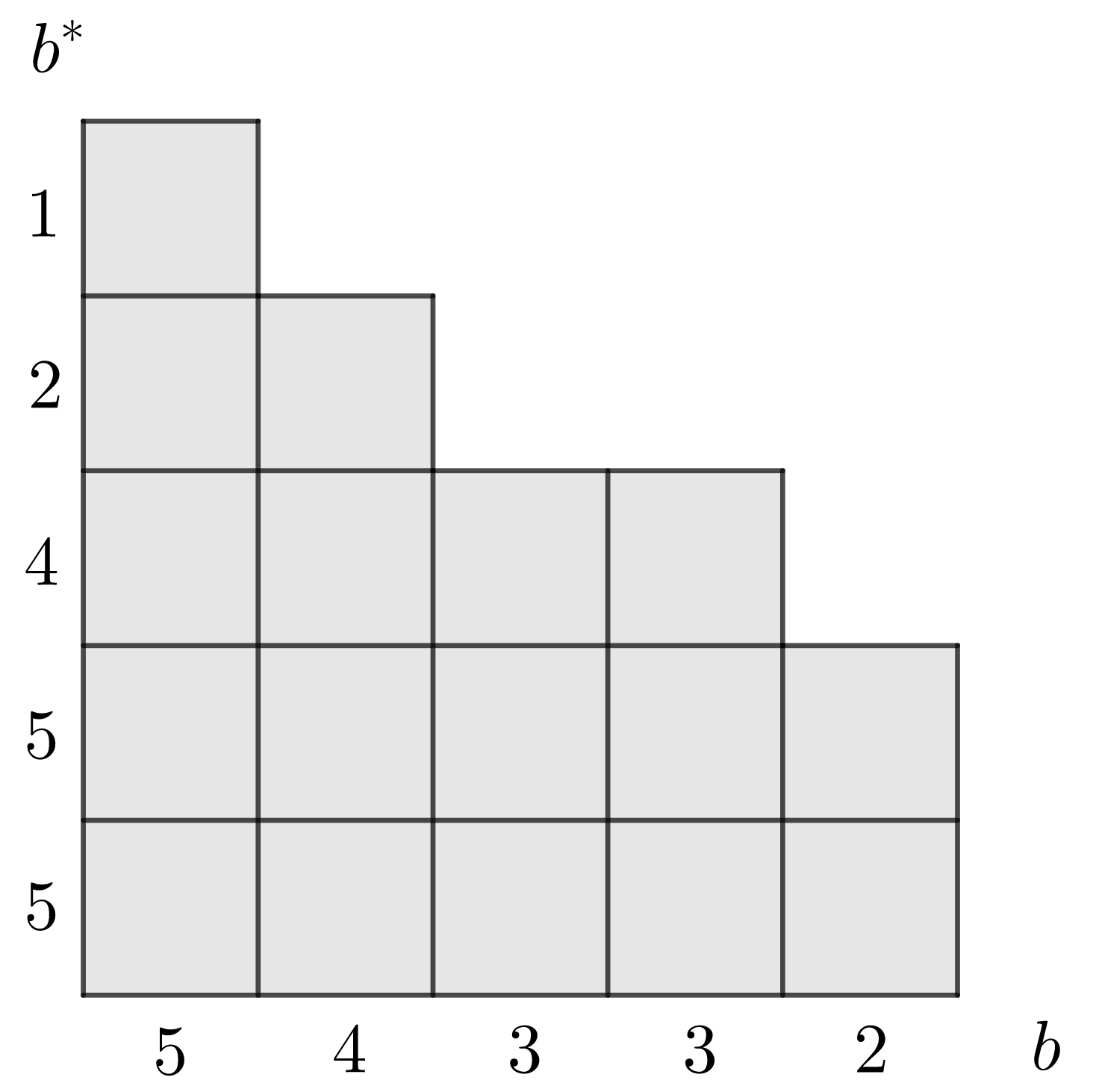}
\caption{ Ferrer diagram  - graphical representation of sequence and its conjugate}
\label{fig:Ferrer}  \end{figure}  

\begin{defi}For real sequneces  $x=(x_{i})_{i=1}^{n}$, $y=(y_{j})_{j=1}^{n}$,  we say that   $x$  majorizes  $y$,   and write  $x \succ y$,  if 
$$x_{\pi(1)} \ \ge \ y_{\sigma(1)},$$
$$x_{\pi(1)} + x_{\pi(2)} \ \ge \ y_{\sigma(1)} + y_{\sigma(2)},$$
$$\cdots$$
$$x_{\pi(1)}+x_{\pi(2)}+\cdots + x_{\pi(n-1)} \ \ge \ y_{\sigma(1)} + y_{\sigma(2)} + \cdots + y_{\sigma(n-1)},$$
$$x_{\pi(1)}+x_{\pi(2)}+\cdots + x_{\pi(n)} \ = \ y_{\sigma(1)} + y_{\sigma(2)} + \cdots + y_{\sigma(n)},$$
where  $\pi$   and  $\sigma$  are such permutations of  $\{1,2,\cdots,n\}$,   that 
$$x_{\pi(1)} \ \ge \ x_{\pi(2)} \ \ge  \ \cdots \ \ge \ x_{\pi(n)}, $$
$$y_{\sigma(1)} \ \ge \ y_{\sigma(2)} \ \ge \cdots \ \ge \ y_{\sigma(n)}.$$

\end{defi}

\begin{theo} \label{Ross-Geller} A pair  $(a,b)$, where  $n\in \mathbb{Z}_{+}$,   $a=(a_{i})_{i=1}^{n}$, $b=(b_{j})_{j=1}^{n}$,   and  
$$n  \ \ge  \  a_{1}, \  a_{2},  \ \cdots , \ a_{n} \  \ge  \ 0,$$   
$$n  \ \ge \  b_{1}, \  b_{2}, \   \cdots , \  b_{n} \  \ge \  0,$$   
is bigraphic, if and only if   $b^{*}  \succ  a$.
\end{theo}
\noindent The next lemma is well known as Karamata's or majorization inequality, see for instance  \cite{wędrówki}.

\begin{lemma}\label{Karamata} Let  $f:\mathbb{R}\rightarrow \mathbb{R}$  be a convex function and assume that two real sequences  $x=(x_{i})_{i=1}^{n}$,  $y=(y_{j})_{j=1}^{n}$  satisfy  $x\succ y$.  Then we have  
$$\sum_{i=1}^{n}f(x_{i}) \ \ge \ \sum_{j=1}^{n}f(y_{j}).$$   \end{lemma}

\begin{theo} \label{b*} For all  $k  \in  \mathbb{Z}_{+}$,  $k\ge 3$,  we have
$$\sup_{n\in\mathbb{Z}_{+}} \ \sup_{\mathcal{B}(n,n)} \ \frac{1}{n^{2+k}} \sum_{i,j=1}^{n} |\deg(x_{i})-\deg(y_{j})|^{k} \ = \ \sup_{n\in \mathbb{Z}_{+}} \  \ \sup_{(b_{j})_{j=1}^{n}  \in  \{0,1,\cdots,n \}^{n}} \ \frac{1}{n^{2+k}} \sum_{i,j=1}^{n}|b_{i}^{*}-b_{j}|^{k}.$$
 \end{theo}
\noindent \textit{Proof}: From  Theorem \ref{Ross-Geller}  we know that
$$\sup_{n\in\mathbb{Z}_{+}} \ \sup_{\mathcal{B}(n,n)} \ \frac{1}{n^{2+k}} \sum_{i,j=1}^{n} |\deg(x_{i})-\deg(y_{j})|^{k} \ = \ \sup_{n\in\mathbb{Z}_{+}} \ \sup_{\substack{ (a_{i})_{i=1}^{n}\in\{0,1,\cdots,n\}^{n} \\ (b_{j})_{j=1}^{n}\in\{0,1,\cdots,n\}^{n} \\ b^{*} \ \succ \ a}} \ \frac{1}{n^{2+k}} 
\sum_{i,j=1}^{n}|a_{i}-b_{j}|^{k}.$$
Let us fix  $n, a, b$  for the time being. After rearrangement, we get 
$$\sum_{i,j=1}^{n}|a_{i}-b_{j}|^{k} \ = \ \sum_{j=1}^{n}\sum_{i=1}^{n}|\{p:b_{p}=j\}|\cdot|a_{i}-j|^k \ = \ \sum_{j=1}^{n}(b_{j}^{*}-b_{j+1}^{*})\sum_{i=1}^{n}|a_{i}-j|^{k},$$
where we put  $b^{*}_{n+1}=0$   for convenience. Note that  $f_{j}(x)=|x-j|^{k}$  is a convex function in  $x$  for all $k\ge3$  and  $b^{*} \succ a$  from the assumption. By  Lemma \ref{Karamata},  for all \ $j\in\{1,2,\cdots,n\}$,  we have
$$ \sum_{i=1}^{n}|a_{i}-j|^{k} \ \le  \ \sum_{i=1}^{n}|b_{i}^{*}-j|^{k}.$$
Summation over  $j$  yields
$$\sum_{i,j=1}^{n}|a_{i}-b_{j}|^{k}  \ \le  \ \sum_{i,j=1}^{n}|b_{i}^{*}-b_{j}|^{k},$$
and hence we obtain
$$\sup_{n\in\mathbb{Z}_{+}} \ \sup_{\mathcal{B}(n,n)} \ \frac{1}{n^{2+k}} \sum_{i,j=1}^{n} |\deg(x_{i})-\deg(y_{j})|^{k} \ \le \ \sup_{n\in \mathbb{Z}_{+}} \  \ \sup_{(b_{j})_{j=1}^{n}  \in  \{0,1,\cdots,n \}^{n}} \ \frac{1}{n^{2+k}} \sum_{i,j=1}^{n}|b_{i}^{*}-b_{j}|^{k}.$$
To prove the opposite inequality,  it is enough to verify that for every  $n\in \mathbb{Z}_{+}$  and every
$$(b_{j})_{j=1}^{n}  \in  \{0,1,\dots,n \}^{n},$$
the pair  $(b^{*},b)$  is bigraphic.  From  Theorem \ref{Ross-Geller},  this is equivalent to
$$b^{*}  \succ \ b^{*},$$
which is clearly true.  \ \ \ $\square$

\subsection{Ferrer diagrams and an upper bound}

Let us start by introducing the following notation:   for  $1 \ge \delta \ge 0$,
$$\epsilon(\delta)=\sup_{(X,Y)\in \mathcal{C}} \mathbb{P}(|X-Y|\ge \delta).$$
In \cite{filo},   \textit{Theorem 18.1},   it was proved that for all  $\delta \in [0,1]$
\begin{equation}
\epsilon(\delta) \ \le \ \big[2(1-\delta)\big]\wedge 1.  \label{eps} 
\end{equation}
By the Fubini's theorem we have the so-called "layer-cake" represantation
\begin{equation}
 \mathbb{E}|X-Y|^{k} =  \int_{0}^{1} ku^{k-1}\cdot \mathbb{P}(|X-Y|\ge u) \ \mathrm{d}u. \label{l-c}  
\end{equation}
Using  (\ref{eps})  and  (\ref{l-c}), we see that for all  $(X,Y)\in \mathcal{C}$  and  $k>0$

\begin{equation}
  \mathbb{E}|X-Y|^{k} \  \le  \   \int_{0}^{1} ku^{k-1}\cdot \epsilon(u) \ \mathrm{d}u   \ =  \ \frac{2-2^{-k}}{1+k}. \label{upbo} 
\end{equation}
The upper bound  (\ref{upbo})  has been considered already by  \textit{Burdzy}  and  \textit{Pitman}  in  \cite{pitman}.   In this section we will reprove this result with additional assumption of independence. Hence, our result is weaker, but the approach we take is different. The reader should treat this section as a soft introduction to the combinatrial ideas that will be studied further in the next chapter. 

\begin{defi}  \label{FD} For   $n\in \mathbb{Z}_{+}$  and any integer sequnece   $b=(b_{i})_{i=1}^{n}$,  with
$$n  \ \ge  \  b_{1} \ \ge \  b_{2} \ \ge  \ \cdots \ \ge \ b_{n} \  \ge  \ 0,$$ 
we define the corresponding   Ferrer diagram    as the   $n\times n$   binary  matrix,  such that
\begin{itemize}
\item its column sums,   starting from the left,  are  $b_{1}, \cdots, b_{n},$  respectively,
\item for every fixed column, all ones  are below  all zeros.
\end{itemize}
\end{defi}
\noindent Rather then thinking in terms of  $0-1$  matrices, we will visualise Ferrer diagrams as square grids with empty or filled cells. 
For example, the Ferrer diagram of the sequence  $b=(5,4,3,3,2)$  is ilustrated in  \textit{Figure }~\ref{fig:Ferrer}.  Note, that the conjugate sequence  $b^{*}$  can now be easily interpreted as the row sums of the Ferrer diagram of  $b$.  In the example given above, we have  $b^{*}=(5,5,4,2,1)$.

\begin{theo} \label{kreska} For all  $k \in \mathbb{Z}_{+} $,  $k\ge2$   we have
$$\sup_{(X,Y)\in \mathcal{C}_{\mathcal{I}}}\mathbb{E}|X-Y|^{k}  \ \le  \ \frac{2-2^{-k}}{1+k}.$$ \end{theo}
\noindent \textit{Proof}: Fix any  $n\in \mathbb{Z}_{+}$  and a decreasing, integer  sequence   $b=(b_{i})_{i=1}^{n}$,  just as in Definition
\ref{FD}. By  Theorems  \ref{b*}  and  \ref{Bk},  it is enough to check that
$$\frac{1}{n^{2+k}} \sum_{i,j=1}^{n}|b_{i}^{*}-b_{j}|^{k} \ = \ \frac{1}{n^{2}} \sum_{i,j=1}^{n}\Bigg|\frac{b_{i}^{*}-b_{j}}{n}\Bigg|^{k} 
\  \le \  \frac{2-2^{-k}}{1+k}. $$
We will start by constructing  $(X,Y)\in\mathcal{C}_{\mathcal{I}}$,  such that 
\begin{equation} \mathbb{E}|X-Y|^{k} \ = \ \frac{1}{n^{2}} \sum_{i,j=1}^{n}\Bigg|\frac{b_{i}^{*}-b_{j}}{n}\Bigg|^{k}. \label{krowa} \end{equation}
For this purpose, take  the Ferrer diagram of  $b$,  and rescale it so that it is  contained it in the unit square
$$I \ = \ \Big\{(u,v): u,v\in [0,1]\Big\}. $$
We shall, from now on, work with the probability space 
$$\Big(I, \mathcal{L}_{ \lambda \otimes \lambda }(I), \lambda \otimes \lambda \Big),$$
 where  $\lambda$  is simply one-dimensional  lebesgue measure on  $[0,1]$,  $\lambda \otimes \lambda$  stands for the product measure on  $I$ and   $\mathcal{L}_{ \lambda \otimes \lambda }(I)$   is an apropriate product $\sigma$-field. Set
$$U(u,v) =  u,$$
$$V(u,v) = v.$$
Therefore   $U, V    \sim    \mathcal{U}[0,1]$,    and   $U\perp V$.   Moreover, by   $A$  let us  denote the staircase-shaped region obtained by uniting  all of the filled cells in the rescaled diagram. We can now define
$$X(u,v)  \ \equiv \  \textbf{x}(U(u,v)) \ := \ \lambda \Bigg(\Big(A\setminus \mathrm{bd}(A)\Big)\cap \Big(u\times[0,1]\Big)\Bigg),$$
$$Y(u,v) \ \equiv \ \textbf{y}(V(u,v)) \ := \ \lambda \Bigg(\Big(A\setminus \mathrm{bd}(A)\Big) \cap \Big([0,1]\times v\Big)\Bigg),$$
where   $\mathrm{bd}(A)$  stands for boundary of  $A$.  This notation indicates that the random variables $X$  and  $Y$  can be also treated as a deterministic (borel) functions 
$$\textbf{x}, \  \textbf{y} \ : \ [0,1] \ \longrightarrow \ [0,1],$$
of random  $U$ \ and  $V$.  Clearly  $X\perp Y$  and the condition  (\ref{krowa})  holds.  Finally, one can also check that
$$X(u,v) \ =  \ \mathbb{E}(\mathds{1}_{A\setminus \mathrm{bd}(A)}|U=u),$$
$$Y(u,v)  \ = \  \mathbb{E}(\mathds{1}_{A\setminus \mathrm{bd}(A)}|V=v),$$
but we postpone the formal verification till the next section.
\begin{figure}[H]
\centering
\includegraphics[width=66mm]{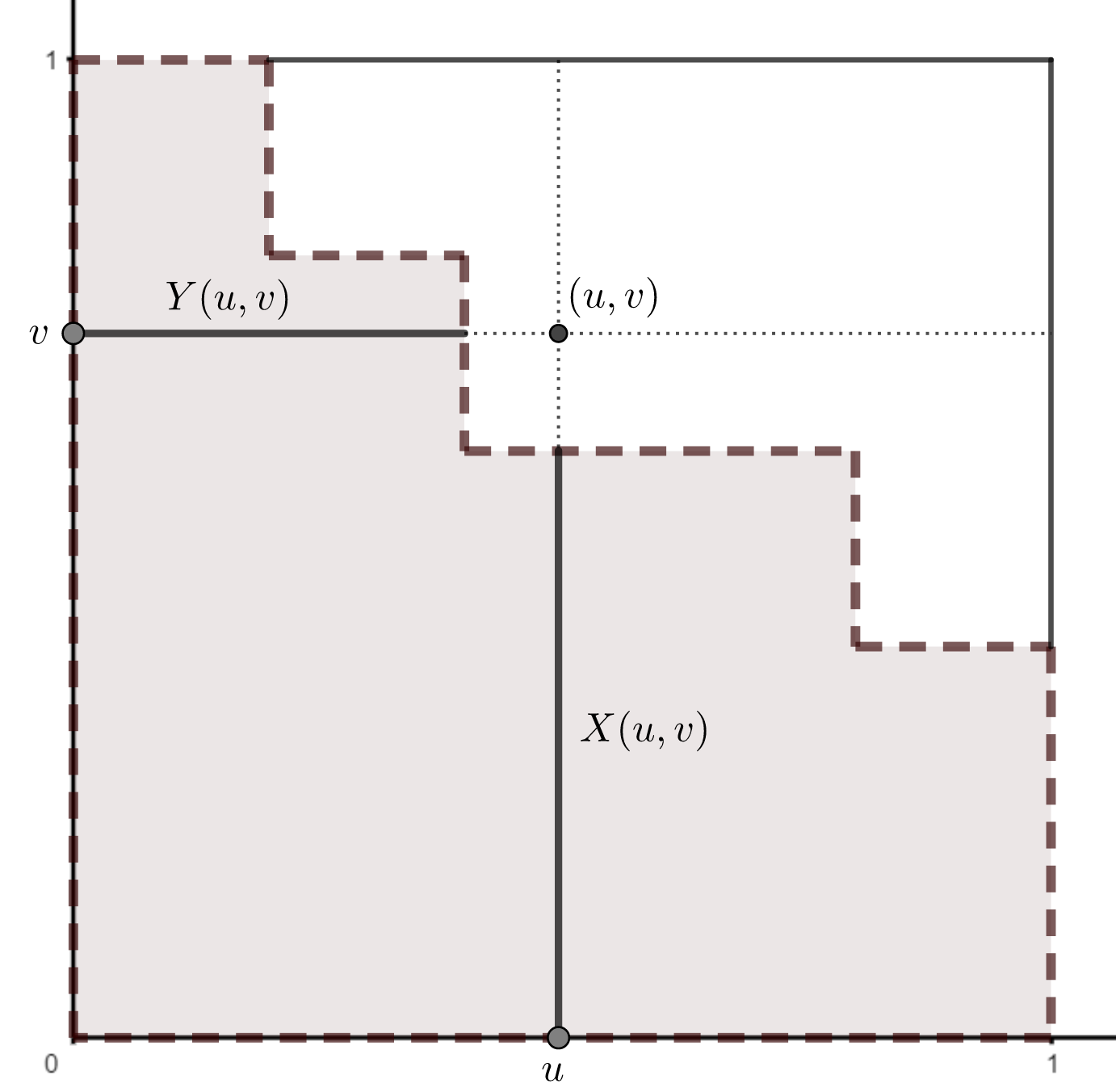}
\caption{ Ferrer diagram    and associated  $(X,Y)\in \mathcal{C}_{\mathcal{I}}$}
\label{fig:Ferrer2}  \end{figure}

\begin{figure}[H]
\centering
\includegraphics[width=125mm]{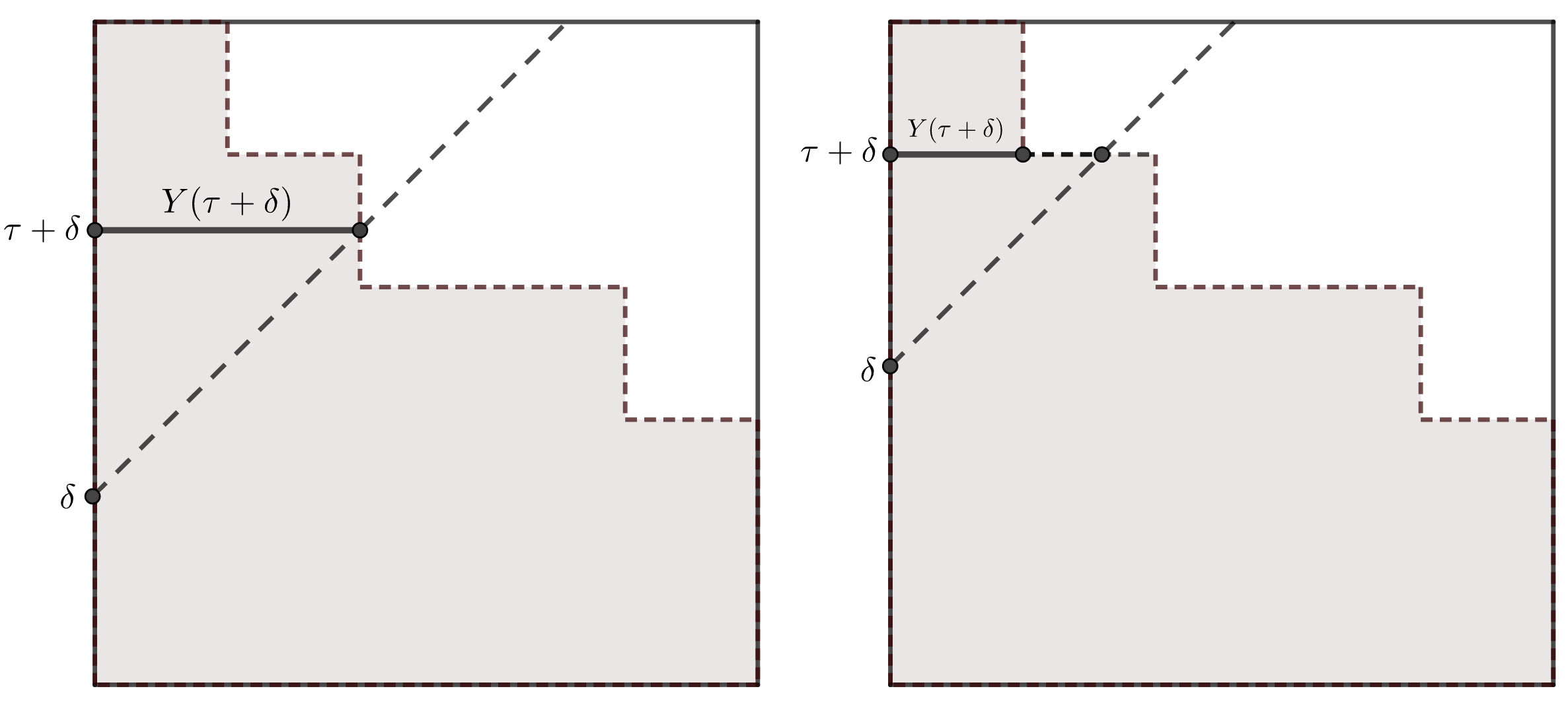}
\caption{ graph of $f$ intersects either horizontal or vertical part of the "staircase"}
\label{fig:intersection}  \end{figure}

\noindent
In this setting we will show the  (\ref{eps})   inequality, namely
$$\mathbb{P}(|X-Y|> \delta ) \ \le \  2(1-\delta),$$
for all  $\delta\in [0,1]$.  This will also establish  (\ref{upbo})  and conclude the proof.  We start by writing
$$\mathbb{P}(|X-Y|> \delta) \  =   \ \mathbb{P}[X> (Y+\delta) ] \ + \ \mathbb{P}[Y> (X+\delta)].$$
By the symmetry of the problem, it is sufficient to demonstrate that
$$\mathbb{P}[X> (Y+\delta)] \ \le \ \ 1-\delta.$$
For any  $\tau \in [0,1]$,  considering the intersection with  $\{Y < \tau \}$,  gives
$$\mathbb{P}[X > (Y+\delta) ] \  = \  \mathbb{P}[X > ( Y+\delta),  \ Y < \tau] \ + \ \mathbb{P}[X >  (Y+\delta), \ Y \ge \tau ] $$
$$\le \ \ \mathbb{P}[Y < \tau] \ + \ \mathbb{P}[X> (\tau + \delta)].$$
Let us now consider the following linear function of  $v$
$$f(v)=v-\delta$$
and let  $\tau$  be such,   that  $\tau+\delta$  is an argument at which the graph of   $f$  and the "staircase" part of  $\mathrm{bd}(A)$  intersect  -  see Figure ~\ref{fig:intersection}. 

\noindent
Clearly,  we have two possible scenarios.  Firstly, graph of $f$ may intersect vertical part of the boundary. In this case, let us note that 
$$\textbf{y}(\tau+\delta) \ = \ f(\tau+\delta) \ = \ (\tau + \delta)-\delta \ = \ \tau.$$
On the other hand, assume that $f$ intersects horizontal part of the "staircase". Having this in mind, recall that we have omitted  $\mathrm{bd}(A)$  in the construction of  $(X,Y)$.  This  leaves us with
$$\textbf{y}(\tau+\delta) \ = \   \lambda \Bigg(\Big(A\setminus \mathrm{bd}(A)\Big) \cap \Big([0,1]\times (\tau+\delta)\Big)\Bigg) \ \le \ f(\tau+\delta) \ = \ \tau.$$
Thus, either way, we get   $\textbf{y}(\tau+\delta)  \le  \tau$.  Again, from construction of   $(X,Y)$   and omission of boundary, we have
$$\mathbb{P}[X> (\tau + \delta)]  \ =  \ \textbf{y}(\tau+\delta) \ \le \  \tau.$$
Hence, it remains to check, that
$$\mathbb{P}[Y< \tau] \ \le  \ 1-\tau-\delta ,$$
or equivalently, that
$$\mathbb{P}[Y\ge \tau]  \ \ge  \ \tau + \delta.$$
 Luckily, observe that
$$\mathbb{P}[Y\ge \tau] \  = \ \mathbb{P}[\textbf{y}(V)\ge \tau],$$
and   $\textbf{y}(v)  \ge  \tau$ for all  $v  < \tau+\delta$. \ \ \  $\square$


\section{Upper bounds and a novel approach}

We  will continue exploiting combinatorial nature of  Ferrer diagrams. We shall start by introducing more flexible definitions. 

\subsection{Generalization of Ferrer diagrams}

\begin{defi} \label{Ff} From now on, by  (generalised) Ferrer diagram,  we shall mean a set
$$F_{f} \  = \  \{(u,v)\in[0,1]^2 \ :  \ v < f(u) \} ,$$
where
$$f:[0,1]\rightarrow [0,1],$$
is  any  weakly decreasing step function that takes finitely  many differernt values;  let us denote the set of such step functions as  
 $\mathbb{STEP}.$ \end{defi}
\noindent In the next definition, we formalise an idea used in the proof of  Theorem \ref{kreska}  -  compare Figure ~\ref{fig:Ferrer2}.

\begin{defi} \label{FfXY}For any  Ferrer diagram  $F_{f}$  we define a pair of associated random variables
$$ (X_{f},Y_{f}) \in \mathcal{C}_{\mathcal{I}} \ \ \ \ \text{defined on a probability space}  \ \ \ \Big(I, \mathcal{L}_{ \lambda \otimes \lambda }(I), \lambda \otimes \lambda \Big),$$
where  $I  =  \Big\{(u,v): u,v\in [0,1]\Big\}$, \ by
\begin{equation}X_{f}(u,v)  \ \equiv \  \textbf{x}_{f}(U(u,v)) \ := \ \lambda \Bigg(\Big(F_f \setminus \mathrm{bd}(F_f)\Big)\cap \Big(u\times[0,1]\Big)\Bigg), \label{defixf} \end{equation}
\begin{equation}Y_{f}(u,v) \ \equiv \ \textbf{y}_{f}(V(u,v)) \ := \ \lambda \Bigg(\Big(F_f\setminus \mathrm{bd}(F_f)\Big) \cap \Big([0,1]\times v\Big)\Bigg), \label{defiyf} \end{equation}
where
$$U(u,v)  =  u,$$
$$V(u,v) =  v.$$
\end{defi}
\noindent We shall  prove that  $(X_f,Y_f)$  defined by   (\ref{defixf})  and   (\ref{defiyf})  does satisfy   $(X_f,Y_f)\in \mathcal{C}_{\mathcal{I}}$.
\\ \\ \noindent \textit{Proof}: Clearly   $U, V    \sim    \mathcal{U}[0,1]$    and   $U\perp V$. This gives  $X\perp Y$.  It is therefore enough to check that
\begin{equation}X_{f}(u,v)  \ \equiv \  \textbf{x}_{f}(U(u,v)) \ \  = \ \ \mathbb{E}(\mathds{1}_{F_f\setminus \mathrm{bd}(F_f)}|U=u), \label{tylkoto} \end{equation}
\begin{equation}Y_{f}(u,v)  \ \equiv \  \textbf{y}_{f}(V(u,v))  \ \ = \ \ \mathbb{E}(\mathds{1}_{F_f\setminus \mathrm{bd}(F_f)}|V=v). \end{equation}
This gives  $(X_f, Y_f)\in \mathcal{C}$.  We limit ourselves to showing  (\ref{tylkoto}).   It is  straightforward to check that  $\textbf{x}_{f}$  is a borel function. Thus by  (\ref{defixf})  we have
$$X_f = \textbf{x}_f(U),$$
and   $X_f$ is   $\sigma(U)$   measurable. It remains to verify that
for every   $A \in \sigma(U)$,  we have
$$\int_{A}X_f  \ \mathrm{d}\mathbb{P} \  =  \ \int_{A}\mathds{1}_{F_f\setminus \mathrm{bd}(F_f)}  \ \mathrm{d}\mathbb{P}.$$
The condition  $A\in \sigma(U)$  is equivalent to  $A=\tilde{A}\times [0,1]$,  for some  $\tilde{A}\in \mathcal{L}([0,1])$.   We can write 
 $$\textbf{x}_{f}(u) \ = \ \lambda \Bigg(\Big(F_f \setminus \mathrm{bd}(F_f)\Big)\cap \Big(u\times[0,1]\Big)\Bigg) \ = \ \int_{0}^{1}\mathds{1}_{F_f\setminus \mathrm{bd}(F_f)}(u,v) \ \mathrm{d}\lambda(v),$$
and hence, by Fubini's theorem, we have 
 $$\int_{A}\mathds{1}_{F_f\setminus \mathrm{bd}(F_f)}  \ \mathrm{d}\mathbb{P} \  =  \ \int_{\tilde{A}} \Bigg [\int_{[0,1]} \mathds{1}_{F_f\setminus \mathrm{bd}(F_f)}(u,v) \ \mathrm{d}\lambda(v)\Bigg] \mathrm{d}\lambda(u)$$
$$ = \  \int_{\tilde{A}} \textbf{x}_{f}(u) \ \mathrm{d}\lambda(u) \  =  \ \int_{A} X_{f}(u,v) \ \mathrm{d}\lambda_{2} (u,v) \ \ = \ \ \int_{A} X_f  \ \mathrm{d}  \mathbb{P},$$ 
as required. \ \ \ $\square$ 
\\ \\ \noindent Although it is rather obvious, let us state the following 

\begin{prop}\label{stepować} For all  $k  \in  \mathbb{Z}_{+}$,  $k\ge 3$,  we simply have
$$ \sup_{f \ \in \ \mathbb{STEP}} \ \mathbb{E}|X_{f}-Y_{f}|^{k}  \ =  \ \sup_{(X,Y)\in \mathcal{C}_{\mathcal{I}}}\mathbb{E}|X-Y|^{k}.$$ \end{prop}
\noindent \textit{Proof}: For all  $f  \in   \mathbb{STEP}$, by definition 
$$(X_f, Y_f) \ \in \mathcal{C}_{\mathcal{I}},$$  and thus, the inequality 
$$ \sup_{f \ \in \ \mathbb{STEP} } \ \mathbb{E}|X_{f}-Y_{f}|^{k} \  \le  \ \sup_{(X,Y)\in \mathcal{C}_{\mathcal{I}}}\mathbb{E}|X-Y|^{k},$$
is clear. The opposite inequality follows from the same argument as the proof of  Theorem \ref{kreska}. \ \ \ $\square$ 

\subsection{Sharpening a layer-cake upper bound}
In this section we will continue our analysis of upper bounds generated by the layer-cake representation  (\ref{l-c}). 
Let us start with a brief overview of the relevant results. As already mentioned, in  \cite{filo}   it was proved that for  $\delta \in (\frac{1}{2},1],$  we have
$$\sup_{(X,Y)\in \mathcal{C}} \mathbb{P}(|X-Y|\ge \delta) \  \le  \ 2(1-\delta).$$
This result has been lately greatly  improved by  \textit{Burdzy}  and  \textit{Pal}.  In  \cite{contra}   they proved that for 
 $\delta \in (\frac{1}{2},1]$: 
$$\sup_{(X,Y)\in \mathcal{C}} \mathbb{P}(|X-Y|\ge \delta) \  = \ \frac{2(1-\delta)}{2-\delta}.$$
Moreover, for all  $\delta$  in this range, one can find such pairs   $(X_{\delta},Y_{\delta})\in \mathcal{C}$  for which the equality is attained. It is however important to note that for  $\delta<1$,  those variables turn out to be dependent.
It is relatively easy to check that for  $\delta \in (\frac{1}{2},1]$:
$$\sup_{(X,Y)\in \mathcal{C}_{\mathcal{I}}(2)}\mathbb{P}(|X-Y|\ge \delta) \  =  \ 2\delta(1-\delta).$$
Based on this premise, \textit{Burdzy}  and  \textit{Pitman}  conjectured in  \cite{pitman},  that for  $\delta \in (\frac{1}{2},1]$, we have
\begin{equation}
\sup_{(X,Y)\in \mathcal{C}_{\mathcal{I}}}\mathbb{P}(|X-Y|\ge \delta)  \ =  \ 2\delta(1-\delta).  \label{hipoteza} 
\end{equation}
In this chapter we will prove a result, that is strikingly similar. Namely, for  $\delta \in (\frac{1}{2},1]$, we have
 \begin{equation}
\sup_{f \ \in \ \mathbb{STEP}}\mathbb{P}(|X_f-Y_f|> \delta) \  =  \ 2\delta(1-\delta).  \label{synteza} 
\end{equation}
In this place let us highlight that, thanks to  Proposition  \ref{stepować},  both  (\ref{hipoteza})  and  (\ref{synteza})  generate exactly
the same layer-cake bound on 
$$\sup_{(X,Y)\in \mathcal{C}_{\mathcal{I}}}\mathbb{E}|X-Y|^{k},$$
for all  $k\in \mathbb{Z}_{+},   k\ge3.$  We begin with the simple, but useful observation.

\begin{figure}[H]
\centering
\includegraphics[width=125mm]{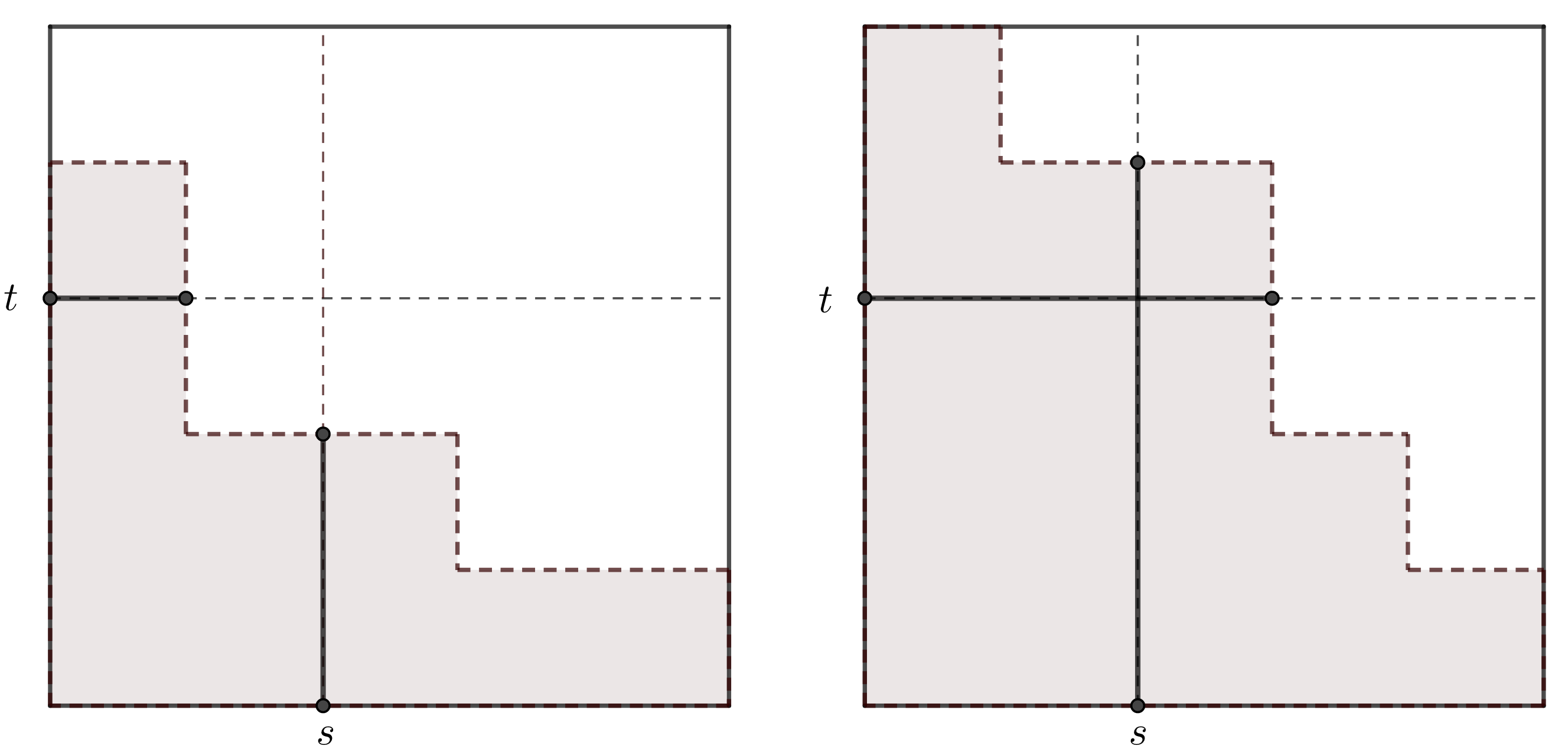}
\caption{  either \ \  $(\textbf{x}_f(s)\le t) \ \wedge \ (\textbf{y}_f(t)\le s)$ \ \   or  \ \  $(\textbf{x}_f(s)\ge t) \  \wedge  \ (\textbf{y}_f(t)\ge s)$.}
\label{fig:xletyles}  \end{figure}

\begin{lemma} \label{tercio} Fix  $f  \in   \mathbb{STEP}$  and consider the associated  $(X_f,Y_f) \in \mathcal{C}_{\mathcal{I}}$.  For  any  $s,t \in [0,1]$,  there are only two possible scenarios:
$$\text{either} \ \ \ \
\left\{ \begin{array}{ll}
\textbf{x}_f(s) \ \le \ t,  \\
\textbf{y}_f(t) \ \le \ s,  
\end{array} \right. 
\ \ \ \ \ \text{or} \ \ \ \ \
\left\{ \begin{array}{ll}
\textbf{x}_f(s) \ \ge \ t,  \\
\textbf{y}_f(t) \ \ge \ s. 
\end{array} \right.
$$ \end{lemma}
\noindent \textit{Proof}: This is a direct consequence of the definitions  -  see Figure ~\ref{fig:xletyles}. \ \ \ $\square$

\begin{theo} \label{newbound} For all  $\delta \in (\frac{1}{2},1]$, we have
$$\sup_{f \ \in \ \mathbb{STEP}}\mathbb{P}(|X_f-Y_f|> \delta) \  =  \ 2\delta(1-\delta).$$ \end{theo}
\noindent Proof of this theorem will be based on two lemmas, which follow below.

\begin{lemma}\label{lemacisko} Fix   $\delta \in (\frac{1}{2},1]$  and consider any   $f  \in  \mathbb{STEP}$,   for which one of conditions:
\begin{equation}
\left\{ \begin{array}{ll}
\textbf{x}_f(\delta)  \ \ge  \ \delta,  \\
\textbf{y}_f(\delta) \  \ge  \ \delta,    
\end{array} \right. \label{a}
\end{equation}
or
\begin{equation}
\left\{ \begin{array}{ll}
\textbf{x}_f(1-\delta)  \ \le   \ 1-\delta,  \\
\textbf{y}_f(1-\delta) \ \le  \ 1-\delta,   
\end{array} \right. \label{b}
\end{equation}
is satisfied. Then we have
$$\mathbb{P}(|X_f-Y_f|> \delta) \  \le  \ 2\delta(1-\delta).$$ \end{lemma}
\noindent \textit{Proof}: Let us first consider point  (\ref{a}).  We start by writing
$$\mathbb{P}(|X_f-Y_f|> \delta)  \ =   \ \mathbb{P}[X_f>(Y_f+\delta)] \ + \ \mathbb{P}[Y_f> (X_f+\delta)].$$
By independence, we can evaluate
$$\mathbb{P}[X_f>(Y_f+\delta)] \  \le \  \mathbb{P}[X_f> \delta]\cdot \mathbb{P}[Y_f\le \delta] \  =  \ \textbf{y}_f(\delta)(1-\textbf{x}_f(\delta)).$$
$$\mathbb{P}[Y_f> (X_f+\delta)] \  \le  \ \mathbb{P}[Y_f> \delta]\cdot \mathbb{P}[X_f\le \delta] \  =  \ \textbf{x}_f(\delta)(1-\textbf{y}_f(\delta)).$$
Summing up, we get
$$\mathbb{P}(|X_f-Y_f|> \delta) \  \ \ \le \ \ \  \textbf{y}_f(\delta)(1-\textbf{x}_f(\delta)) \ + \  \textbf{x}_f(\delta)(1-\textbf{y}_f(\delta))$$
$$ \ \ \ \le \ \ \ \sup_{x,y \ \in \ [\delta, 1]}  g(x,y),$$
where  $g(x,y)  =   y(1-x)  +  x(1-y).$  Since  $\delta>\frac{1}{2},$  we get
$$\frac{\partial g}{\partial x}(x,y)  \ =  \ 1-2y  \ <  \ 0,$$
$$\frac{\partial g}{\partial y}(x,y) \  = \  1-2x \  < \ 0,$$
and hence
$$\mathbb{P}(|X_f-Y_f|> \delta)   \   \le \  g(\delta, \delta)  =  2\delta(1-\delta), $$
which was to be proved. 
\\ \\The proof of point  (\ref{b})  is similar and we will only sketch it. One can evaluate
$$\mathbb{P}(|X_f-Y_f|> \delta) \  \le \  \mathbb{P}[X_f> 1- \delta]\cdot \mathbb{P}[Y_f\le 1- \delta] \ + \ \mathbb{P}[Y_f> 1- \delta]\cdot \mathbb{P}[X_f\le 1- \delta] $$
$$= \  \textbf{y}_f(1-\delta)(1-\textbf{x}_f(1-\delta)) \ + \  \textbf{x}_f(1-\delta)(1-\textbf{y}_f(1-\delta))$$
$$   \le  \ \sup_{x,y \ \in \ [0, 1-\delta]} g(x,y) \  =  \ 2\delta(1-\delta),$$
where the last line is a consequence of  $1-\delta<\frac{1}{2}$.  \ \ \ $\square$

\begin{lemma} \label{okno} For all   $\delta \in (\frac{1}{2},1]$  and any  $f  \in  \mathbb{STEP}$,  if 
$$ 1-\delta  \ \le  \ \textbf{y}_f(\delta)  \ \le  \ \delta,$$ 
then
$$\mathbb{P}(X_f>Y_f+\delta) \  \le \  \delta(1-\delta).$$
The same holds true for  $X_f$  and  $Y_f$   with switched roles.\end{lemma}
\noindent \textit{Proof}: Just as previously, due to independence, we can write
$$\mathbb{P}(X_f>Y_f+\delta) \  \le  \  \mathbb{P}[X_f> \delta]\cdot \mathbb{P}[Y_f< 1- \delta].$$
Firstly, note that
$$\mathbb{P}[X_f>\delta] \  =  \ \textbf{y}_f(\delta) \  \le  \ \delta.$$
Secondly, thanks to monotonicity of   $Y_f$,  we have
$$\textbf{y}_{f}(\omega) \  \ge  \ \textbf{y}_{f}(\delta )  \ \ge  \ 1-\delta \ \ \ \ \ \ \ \text{for all} \ \ \ \ \ \ \ \omega \ \le \ \delta,$$
so  $\{\textbf{y}_f <1-\delta\}  \subset  (\delta,1],$  and hence
$$ \mathbb{P}[Y_f< 1- \delta] \  \le  \ 1-\delta,$$
which completes the proof. \ \ \ $\square$
\\ \\ \textit{Proof  of   Theorem  \ref{newbound}}:   Fix  $\delta \in (\frac{1}{2},1]$  and any  $f\in \mathbb{STEP}$.   It is enough to check, that
$$\mathbb{P}(|X_f-Y_f|> \delta) \  \le \  2\delta(1-\delta).$$
By  Lemmas  \ref{tercio}  and  \ref{lemacisko}  (\ref{a}),  we can assume, that
\begin{equation}
\left\{ \begin{array}{ll}
\textbf{x}_f(\delta) \ \le \ \delta,  \\
\textbf{y}_f(\delta) \ \le \ \delta.  
\end{array} \right. \label{aleph}
\end{equation}
By the  Lemma \ref{tercio}  again,  there are only  $4$  possible scenarios:
\begin{equation}
 \left\{ \begin{array}{ll}
\textbf{x}_f(1-\delta) \ \ge \ \delta,  \\
\textbf{y}_f(\delta) \ \ge \ 1-\delta,  
\end{array} \right. 
\ \ \ \text{and} \ \ \ 
\left\{ \begin{array}{ll}
\textbf{y}_f(1-\delta) \ \ge \ \delta,  \\
\textbf{x}_f(\delta) \ \ge \ 1-\delta.  
\end{array} \right. \label{A}
\end{equation}

\begin{equation}
 \left\{ \begin{array}{ll}
\textbf{x}_f(1-\delta) \ \ge \ \delta,  \\
\textbf{y}_f(\delta) \ \ge \ 1-\delta,  
\end{array} \right. 
\ \ \ \text{and} \ \ \ 
\left\{ \begin{array}{ll}
\textbf{y}_f(1-\delta) \ \le \ \delta,  \\
\textbf{x}_f(\delta) \ \le \ 1-\delta.  
\end{array} \right. \label{B}
\end{equation}

\begin{equation}
 \left\{ \begin{array}{ll}
\textbf{x}_f(1-\delta) \ \le \ \delta,  \\
\textbf{y}_f(\delta) \ \le \ 1-\delta,  
\end{array} \right. 
\ \ \ \text{and} \ \ \ 
\left\{ \begin{array}{ll}
\textbf{y}_f(1-\delta) \ \ge \ \delta,  \\
\textbf{x}_f(\delta) \ \ge \ 1-\delta.  
\end{array} \right. \label{C}
\end{equation}

\begin{equation}
 \left\{ \begin{array}{ll}
\textbf{x}_f(1-\delta) \ \le \ \delta,  \\
\textbf{y}_f(\delta) \ \le \ 1-\delta,  
\end{array} \right. 
\ \ \ \text{and} \ \ \ 
\left\{ \begin{array}{ll}
\textbf{y}_f(1-\delta) \ \le \ \delta,  \\
\textbf{x}_f(\delta) \ \le \ 1-\delta.  
\end{array} \right. \label{D}
\end{equation}
We continue by  inspection, one by one.
\\ \\(\ref{A}).   In view of  (\ref{aleph})  we have
$$ 1-\delta \  \le  \ \textbf{y}_f(\delta)  \ \le  \ \delta,$$ 
$$ 1-\delta \  \le \ \textbf{x}_f(\delta) \  \le  \ \delta.$$ 
By a double use of  Lemma  \ref{okno},   we get
$$\mathbb{P}(|X_f-Y_f|> \delta) \  =   \ \mathbb{P}[X_f>(Y_f+\delta)] \ + \ \mathbb{P}[Y_f> (X_f+\delta)]$$
$$  \le  \ \delta(1-\delta) \ + \  \delta(1-\delta)  \ =  \ 2\delta(1-\delta). \ \ \ \ \triangle $$
\\(\ref{B}).   In view of  (\ref{aleph})  we have
$$ 1-\delta \  \le  \ \textbf{y}_f(\delta)  \ \le  \ \delta.$$
 From  Lemma  \ref{okno},   we have 
$$\mathbb{P}[X_f>Y_f+\delta] \  \le  \ \delta(1-\delta).$$
Furthermore, since to   $\delta \in (\frac{1}{2},1]$,   we can evaluate
$$\mathbb{P}[Y_f>X_f+\delta]  \ \le \ \mathbb{P}(Y_f>\delta)\cdot \mathbb{P}(X_f<\delta).$$
By monotonicity,  we have
$$\textbf{y}_{f}(v) \  \le  \ \textbf{y}_{f}(1-\delta )  \ \le  \ \delta \ \ \ \ \ \ \ \text{for all} \ \ \ \ \ \ \ v \ \ge \ 1- \delta,$$
$$\textbf{x}_{f}(u) \  \ge  \ \textbf{x}_{f}(1-\delta )  \ \ge  \ \delta \ \ \ \ \ \ \ \text{for all} \ \ \ \ \ \ \ u \ \le \ 1-\delta,$$
so 
$$\{\textbf{y}_f >\delta \} \  \subset  \ [0,1-\delta),$$
$$\{\textbf{x}_f <\delta \} \  \subset  \ (1-\delta,1],$$
and hence
$$\mathbb{P}[Y_f>X_f+\delta]  \ \le \  (1-\delta)\delta. \ \ \ \ \triangle$$
\\ \\(\ref{C}). This scenario is analogous to  (\ref{B}).  It is sufficient to change roles of  $X_f$ and  $Y_f$. \ \  $\triangle$
\\ \\(\ref{D}).  Let us start by repeating again the bounds for this scenario. We have
\[
 \left\{ \begin{array}{ll}
\textbf{x}_f(1-\delta) \ \le \ \delta,  \\
\textbf{x}_f(\delta) \ \le \ 1-\delta,  
\end{array} \right. 
\ \ \ \text{and} \ \ \ 
\left\{ \begin{array}{ll}
\textbf{y}_f(1-\delta) \ \le \ \delta,  \\
\textbf{y}_f(\delta) \ \le \ 1-\delta.  
\end{array} \right. \tag{\ref{D}}
\]
At this point,  it is beneficial to graph the constraints given by  (\ref{D}).

\begin{figure}[H]
\centering
\includegraphics[width=65mm]{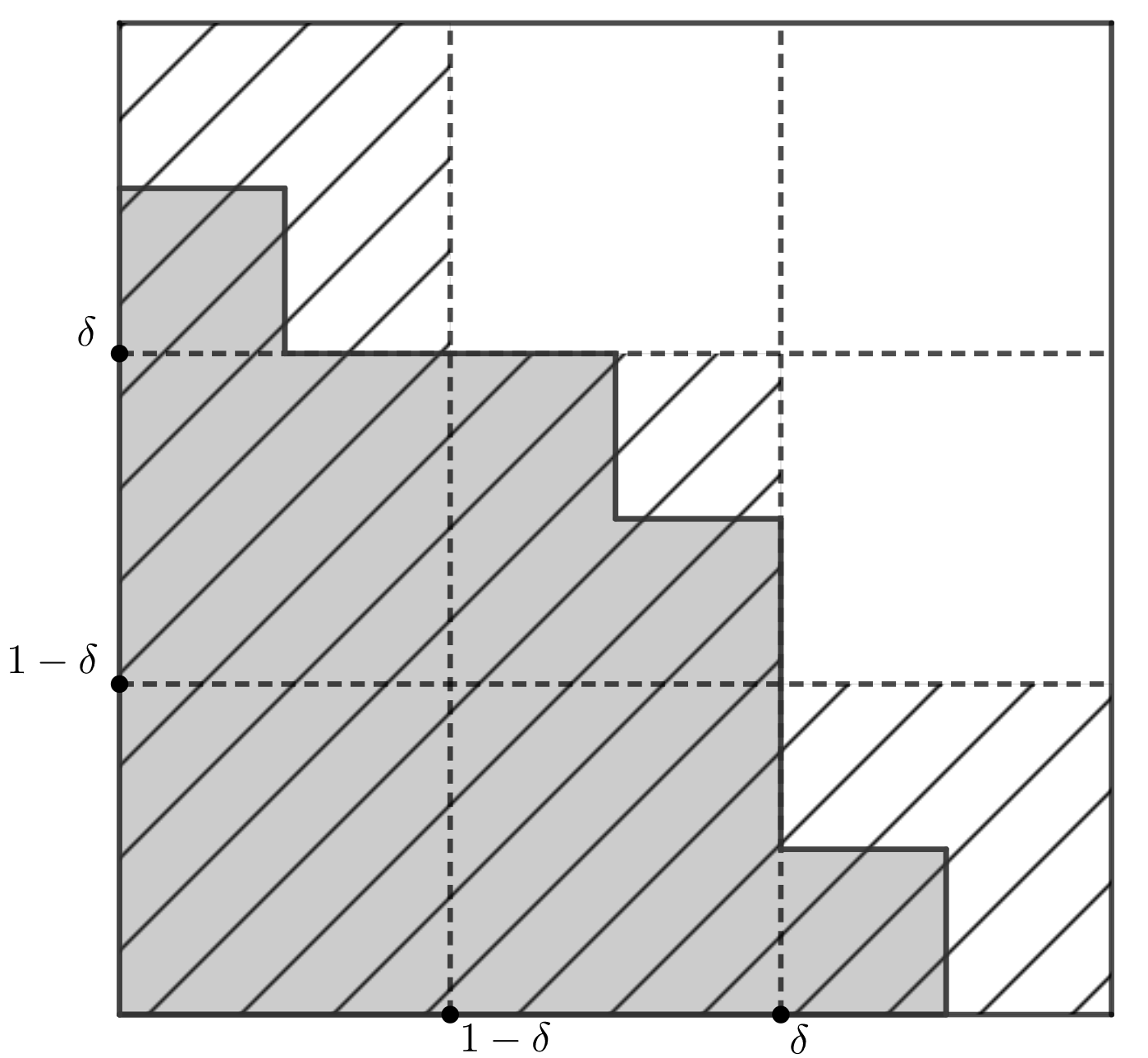}
\caption{ example diagram  $F_f$  meeting constraints given by  (\ref{D}).}
\label{fig:shaded1}  \end{figure} 

\begin{figure}[H]
\centering
\includegraphics[width=65mm]{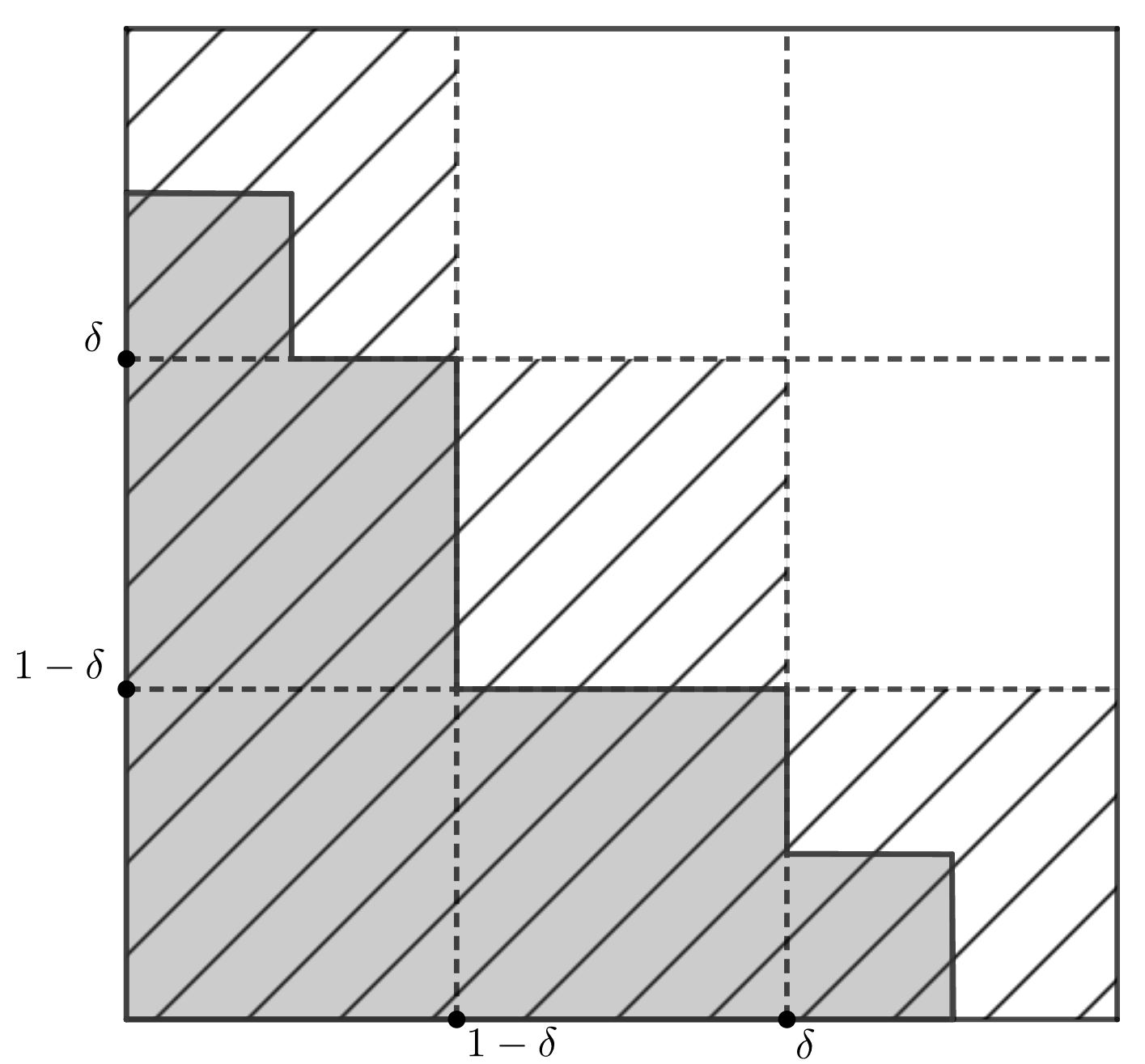}
\caption{diagram $F_{f'}$  is obtained  from   $F_f$  by removing  $[1-\delta, \delta]^2$ .}
\label{fig:shaded2}  \end{figure}

\noindent As in   Figure ~\ref{fig:shaded1}; every diagram  $F_f$  meeting constraints discussed in scenario   (\ref{D}),   must be a subset of hatched region. For any such diagram  $F_f$,  let us now define
$$F_{f'} =  F_{f}  \setminus   [1-\delta, \delta]^2.$$
Put differently, we are removing the  (possibly empty)  interesection  $F_f \cap [1-\delta, \delta]^2$  - see  Figure ~\ref{fig:shaded2}.

\noindent Note that, by construction, the  transformation  $ F_f    \to  F_{f'}$   fulfils both
$$\textbf{x}_{f'} \ \le \  \textbf{x}_f,$$
$$\textbf{y}_{f'} \ \le \ \textbf{y}_{f},$$ 
and
$$\textbf{x}_{f'}(u) \ = \ \textbf{x}_{f}(u) \ \ \ \ \ \ \ \text{for all} \ \ \ \ \ \ \ u \ \in \ \{u: \textbf{x}_{f}(u)>\delta\},$$
$$\textbf{y}_{f'}(v) \ = \ \textbf{y}_{f}(v) \ \ \ \ \ \ \ \text{for all} \ \ \ \ \ \ \ v \ \in \ \{v: \textbf{y}_{f}(v)>\delta\}.$$
Thus, it is straightforward to see, that
$$\mathbb{P}(X_{f}>Y_f+\delta) \ \le \  \mathbb{P}(X_{f'}>Y_{f'}+\delta),$$
$$\mathbb{P}(Y_{f}>X_f+\delta) \  \le  \ \mathbb{P}(Y_{f'}>X_{f'}+\delta),$$
and
$$\mathbb{P}(|X_f-Y_f|>\delta) \  \le  \  \mathbb{P}(|X_{f'}-Y_{f'}|>\delta),$$
as a result.  To complete the proof, it is enough to show, that
$$ \mathbb{P}(|X_{f'}-Y_{f'}|>\delta)  \ \le  \ 2\delta(1-\delta),$$
but this is a direct  consequence of  Lemma   \ref{lemacisko} (\ref{b}). \ \ \ $\triangle$ \  $\square$
\\ \\By  Proposition  \ref{stepować} and  Theorem  \ref{newbound},  we  get the following corollary directly.

\begin{coro} For all  $k  \in  \mathbb{Z}_{+}$,  $k\ge 3$,  we  have
$$\sup_{(X,Y)\in \mathcal{C}_{\mathcal{I}}}\mathbb{E}|X-Y|^k  \ \le \  \int_{0}^{\frac{1}{2}}  kt^{k-1}\mathrm{d}t \ + \  \int_{\frac{1}{2}}^{1} kt^{k-1}\cdot2t(1-t) \mathrm{d}{t},$$
that is
$$\sup_{(X,Y)\in \mathcal{C}_{\mathcal{I}}}\mathbb{E}|X-Y|^k  \ \le   \    2\cdot \frac{k}{(k+1)(k+2)}+ \ 2^{-k} \ - \ 2^{-k-1} \cdot \frac{k(k+3)}{(k+1)(k+2)} .$$\end{coro}

\setlength{\baselineskip}{2ex}

\end{document}